# From Road Congestion to Vehicle-Control Enabled Artificial Traffic Fluids

**Authors:** Iasson Karafyllis[1], Dionysios Theodosis[2,*], Markos Papageorgiou[2,3], Miroslav Krstic[4]

**Affiliations:**
[1]Department of Mathematics, National Technical University of Athens; Athens, 15780, Greece.
[2]Dynamic Systems and Simulation Laboratory, Technical University of Crete; Chania, 73100, Greece.
[3]Faculty of Maritime and Transportation, Ningbo University; Ningbo, China.
[4]Department of Mechanical and Aerospace Engineering, University of California, San Diego, La Jolla, CA 92093-0411, USA.
*Corresponding author. Email: dtheodosis@dssl.tuc.gr

**Abstract:** This article provides an overview of the design of nonlinear feedback Cruise Controllers (CCs) for automated vehicles on lane-free roads. The feedback design problem is particularly challenging because of the various state constraints (e.g., collision-free movement, road geometry, speed limits) as well as the nature of the control objective (globally stabilizing distributed controllers that require measurements from neighboring vehicles only). Therefore, the resulting nonlinear control system is defined on an open set (not necessarily diffeomorphic to a linear space) for which the set of desired equilibria is non-compact. The proposed design of the CCs is based on energy-like control Lyapunov functions which combine potential functions with kinetic energy terms and other appropriate penalty terms. The feedback design in the microscopic level is accompanied by the derivation of the corresponding macroscopic traffic flow models. Explicit relations are established between selectable CC features and the obtained macroscopic traffic flow characteristics. This facilitates the active design of efficient traffic flow with desired properties, i.e., the construction of artificial traffic fluids.

## 1. Introduction

This paper reviews some developments related to a paradigm-changing traffic technology. Methods and tools from many scientific fields are employed, but nonlinear feedback control plays the foremost, enabling role.

The twentieth century gave rise to the advent of the automobile which revolutionized urban and interurban mobility and became the dominant means of transport. During the 1950s, to improve traffic safety for the continuously increasing number of vehicles, parallel traffic lanes were introduced, at the expense of reducing road capacity. But, the commercial success of automobiles increased traffic levels dramatically, with existing road networks unable to serve the increasing demand, leading to omnipresent vehicular traffic congestion that reigns during peak periods on the roads around the globe, [18], [19], [73], causing serious travel delays [26], [74], excessive fuel consumption [23], [29], environmental degradation [47], [58], [107], and decreased traffic safety, [26], [29]. Remarkably, this is in contrast with even simple species which manage their transport problems efficiently, see e.g. [2], [14], [100]. Recently, a wide variety of vehicle automation and communication systems have been developed and partly deployed, which tremendously improve the vehicles' individual capabilities, and this tendency is continuing with the emergence of high-automation vehicles which are validated in real traffic environments, see e.g. [1], [16], [28], [43], [55], [90], [91], [103]. For the not-too-distant future, the automotive industry is preparing for vehicles that communicate with each other and with the infrastructure; and drive automatically, based on own sensors, communications, and appropriate movement control strategies, [6], [17], [45], [67], [71], [77], [87], [90], [101].

In this context, the TrafficFluid concept, a novel paradigm for vehicular traffic at high levels of vehicle automation, was recently proposed in [75], relying on two combined principles: (a) Lane-free traffic, whereby vehicles are not bound to fixed traffic lanes, as in conventional traffic, but may drive anywhere on the 2-D surface of the road; and (b) Vehicle nudging, whereby vehicles influence other vehicles not only behind, but also in front or on the sides of them. These two principles are absent in human driving due to a significant perception limitation: human vision is aligned with one direction at a time, which is naturally the direction of movement, and this has two influential implications:

(i) Traffic lanes were introduced in the 1950s to render human driving easier and safer, as driving on a lane reduces to car-following, i.e., keeping a safe distance to the vehicle in front. On the other hand, lanes must be wide enough to accommodate the widest circulating vehicles, hence they lead to lower lateral occupancy and, as a consequence, to lower flow and capacity. Also, traffic lanes necessitate occasional lane-changing, which is a difficult and accident-prone maneuver, [3], [81].

(ii) Driving is only influenced by traffic conditions downstream, which leads to structurally anisotropic flow and lack of pressure at the macroscopic level.

In summary, traffic lanes and anisotropy, stemming from human perception limitations, are impacting traffic safety and efficiency. Particles in other fluids (liquids or gases) are isotropic and do not form lanes. In living nature (bird flocks [7], fish schools [76]), and even in some human activities (crowd or bicycle-avenue movement, [33], [36]), lanes are not encountered.

Connected Automated Vehicles (CAVs) are different objects, compared to human-driven vehicles, and the road traffic organization may have to evolve to better serve the purpose. In particular, traffic lanes and anisotropy may appear, in a few decades, like fossils from the era of human driving. In the absence of lanes, CAVs may drive smoothly in lateral direction, as abrupt lane-changing displacements become obsolete. Also, CAV driving may be designed to depend not only on vehicles in front, but also on the sides or behind them, like isotropic physical particles.

In this article, we highlight recent progress in 2-D CAV movement strategies and describe the foundations of a new scientific CAV traffic flow, for which we construct the fundamental laws, leveraging notions and tools from diverse scientific fields like Control Theory, Numerical Analysis, Ordinary and Partial Differential Equations, Dynamical Systems, Mathematical Physics, Mechanics of Particles and n-body problems, and Fluid Mechanics. The feedback design problem of CCs for CAVs falls into a particular class of control problems that need to be studied on a set which is not necessarily diffeomorphic to a linear space. Due to various constraints, such as collision avoidance and speed positivity, the system of CAVs under consideration is defined on an open set. Moreover, the fact that CAVs, under decentralized CCs, do not interact with each other when the distance between them is large, adds the additional characteristic that the set of equilibrium points may not be compact. Both those features are rarely encountered in control problems. Here, we discuss how 2-D Cruise Controllers (CCs) for CAVs can be designed via the algorithmic creation of artificial forces that induce the vehicle motion in lane-free traffic with vehicle nudging, rigorously guaranteeing a variety of desired features, including safety in terms of collision avoidance and road boundary respect.

In this context, the resulting controlled CAVs resemble isotropic particles whose movement determines the characteristics of the emerging traffic fluid. Indeed, in the physical world, particle movement determines the characteristics of the emerging fluid (e.g., water, gas etc.), [69]. Notice however that, in the case of CAV traffic flow, the particle movement is largely in our hands to design, e.g., via algorithmic creation of artificial forces that induce the vehicle motion. Thus, under the TrafficFluid concept, we face the question: Which (microscopically feasible, convenient and safe) CAV movement strategy inherits opportune properties to the emerging traffic flow? In other words, instead of modeling traffic flow as a physical process, as done for decades [27], we have the opportunity to actively design it, via appropriate vehicle movement strategies, which implies that we may engineer traffic flow as an efficient *artificial fluid,* see also Fig. 1).

**Notation**

* By $|x|$ we denote both the Euclidean norm of a vector $x \in \mathbb{R}^n$ and the absolute value of a scalar $x \in \mathbb{R}$.

* $\mathbb{R}_+ := [0, +\infty)$ denotes the set of non-negative real numbers.

* By $x'$ we denote the transpose of a vector $x \in \mathbb{R}^n$.

* Let $A \subseteq \mathbb{R}^n$ be an open set. By $C^0(A; \Omega)$, we denote the class of continuous functions on $A \subseteq \mathbb{R}^n$, which take values in $\Omega \subseteq \mathbb{R}^m$. By $C^k(A; \Omega)$, where $k \geq 1$ is an integer, we denote the class of functions on $A \subseteq \mathbb{R}^n$ with continuous derivatives of order $k$, which take values in $\Omega \subseteq \mathbb{R}^m$. When $\Omega = \mathbb{R}$ the we write $C^0(A)$ or $C^k(A)$.

* By $K$ we denote the class of increasing $C^0$ functions $a : \mathbb{R}_+ \to \mathbb{R}_+$ with $a(0) = 0$. By $K_\infty$ we denote the class of increasing $C^0$ functions $a : \mathbb{R}_+ \to \mathbb{R}_+$ with $a(0) = 0$ and $\lim_{s \to +\infty} a(s) = +\infty$. By $KL$ we denote the set of all continuous functions $\sigma : \mathbb{R}_+ \times \mathbb{R}_+ \to \mathbb{R}_+$ with the properties: (i) for each $t \geq 0$ the mapping $\sigma(\cdot, t)$ is of class $K$; (ii) for each $s \geq 0$, the mapping $\sigma(s, \cdot)$ is non-increasing with $\lim_{t \to +\infty} \sigma(s, t) = 0$.

* Let $I \subseteq \mathbb{R}$ be a given interval. $L^\infty(I)$ denotes the set of equivalence classes of measurable functions $f : I \to \Re$ for which $\|f\|_\infty = ess \sup_{x \in I}(|f(x)|) < +\infty$. By $W^{k,\infty}(I)$, where $k \geq 1$ is an integer, we denote the Sobolev spaces of functions $f \in L^\infty(I)$ which have weak derivatives of order $\leq k$, all of which belong to $L^\infty(I)$.

* Let $u : \mathbb{R}_+ \times \mathbb{R} \to \mathbb{R}$, $(t, x) \to u(t, x)$ be any function differentiable with respect to its arguments. We use the notation $u_t(t, x) = \frac{\partial u}{\partial t}(t, x)$ and $u_x(t, x) = \frac{\partial u}{\partial x}(t, x)$ for the partial derivatives of $u$ with respect to $t$ and $x$, respectively. We use the notation $u[t]$ to denote the profile at certain $t \geq 0$, $(u[t])[x] := u(t, x)$, for all $x \in \mathbb{R}$.

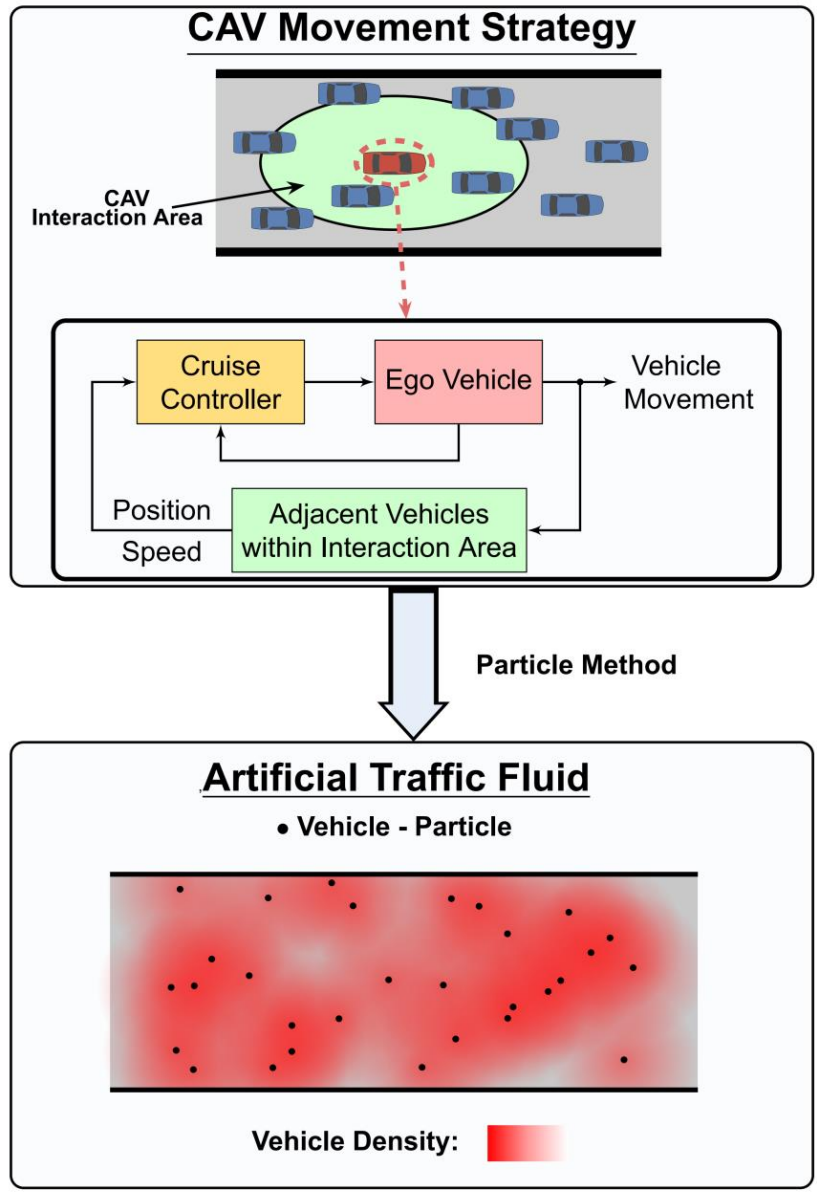

**Fig. 1: Illustration of cruise controller feedback and transition to macroscopic traffic flow models. Top:** The CAV movement strategy involves a feedback loop, where local measurements (distance and speed of adjacent vehicles) from the green vehicle-interaction area are fed to the cruise controller, which, on its turn, determines each vehicle movement. **Bottom:** Via the Particle Method, the emerging CAV traffic flow is derived, with characteristics depending directly on selectable cruise controller structure and parameters, enabling creation of an engineered, hence artificial traffic fluid

## 2. Microscopic Traffic Models

### 2.1 Vehicle Model

The first step of the CC design is the appropriate description of the vehicle dynamics. The bicycle kinematic model is a well-established model for vehicle dynamics since it captures the non-holonomic constraints of the actual vehicle, (see [78], [82]). The equations of motion for each vehicle $i$ are given by the following ODEs (Ordinary Differential Equations)

$$\begin{aligned} \dot{x}_i &= v_i \cos(\theta_i) \\ \dot{y}_i &= v_i \sin(\theta_i) \\ \dot{\theta}_i &= \sigma_i^{-1} v_i \tan(\delta_i) \\ \dot{v}_i &= F_i \end{aligned} \quad (1)$$

for $i = 1,...,n$, where (see Fig. 2) $\sigma_i > 0$ is the vehicle length; $(x_i, y_i)$ is the vehicle's reference point and is placed at the midpoint of the rear axle, with $x_i$ being the longitudinal position and $y_i$ being the lateral position of the vehicle in an inertial frame with Cartesian coordinates $(X,Y)$; $v_i > 0$ is the vehicle speed at the point $(x_i, y_i)$, and $\theta_i \in \left(-\frac{\pi}{2}, \frac{\pi}{2}\right)$ is the vehicle orientation with respect to the $X$ axis. Finally, $\delta_i$ and $F_i$ are the control inputs of the $i$-th vehicle, with $\delta_i$ being the steering angle of the front wheels relative to the orientation $\theta_i$, and $F_i$ being the acceleration.

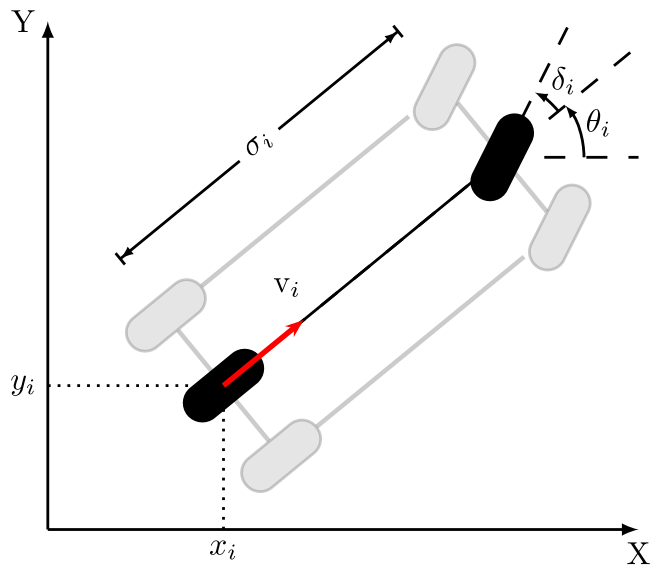

**Fig. 2:** The bicycle kinematic model

## 2.2 Controller Design for Lane-Free Constant-Width Roads

At the microscopic level, each vehicle is driven by a 2-D movement strategy that we call cruise controller (CC). The CC design is a feedback controller design problem. The vehicle's own state, as well as the current distances and speeds of other vehicles and the distances from the road boundaries are measured and "fed back" to the CC, which in turn determines the appropriate vehicle acceleration $F_i$ and steering angle $\delta_i$ (control inputs) that determine the vehicle movement in real-time. Since the range of sensor measurements (or vehicle-to-vehicle communication) is limited, it is important that the required real-time information pertains only to a neighborhood of each vehicle.

The movement of a vehicle is influenced by a plethora of particular and changing characteristics that refer to the vehicle itself and the road infrastructure, such as straight-road of constant or non-constant width, ring-road, urban intersections, on-ramps, off-ramps and more. Individual vehicles must decide on their movement on the basis of real-time information regarding, beyond their own current state, other vehicles and the infrastructure and guarantee that all vehicles: respect the road speed limit $v_{\max} > 0$ and never drive backwards (i.e., $v_i(t) \in (0, v_{\max})$); remain within the road boundaries (i.e., $(x_i(t), y_i(t)) \in \mathbb{R} \times (-a, a)$, where $2a$ is the width of the road); and not turn perpendicular to the road (i.e., $\theta_i(t) \in (-\varphi, \varphi)$ where $\varphi \in \left(0, \frac{\pi}{2}\right]$ is a given constant angle). Since safety is a primary objective for CAVs, vehicles must avoid collision with adjacent vehicles. Defining the distance between two vehicles $i$ and $j$ by

$$d_{i,j} := \sqrt{(x_i - x_j)^2 + p_{i,j}(y_i - y_j)^2} \text{ , for } i, j = 1, \ldots, n \tag{2}$$

with $p_{i,j} > 0$ being constants that satisfy $p_{i,j} = p_{j,i}$ for all $i, j = 1, \ldots, n$, the collision-avoidance objective requires that $d_{i,j}(t) > L_{i,j}$ for all $t \geq 0$ and $i, j = 1, \ldots, n$, $j \neq i$, where $L_{i,j} > 0$ is a safety distance that prevents collisions and satisfies $L_{i,j} = L_{j,i}$ for $i, j = 1, \ldots, n$, $i \neq j$. Note that the distance is defined in such a way that iso-distance curves are ellipses, so that, given the rectangular shape of vehicles, the road width is better exploited in lane-free traffic, see Fig. 3. Finally, each vehicle should use a limited amount of real-time measurements regarding other vehicles and the road boundaries, namely only from a neighborhood of the given vehicle.

The mentioned constraints imply that we must consider the system described by (1) on the open set $\Omega \subset \mathbb{R}^{4n}$ defined by:

$$\Omega := \left\{ w \in \mathbb{R}^{4n} : \begin{array}{l} x_i \in \mathbb{R}, |y_i| < a, i = 1, \ldots, n \\ v_i \in (0, v_{\max}), |\theta_i| < \varphi, i = 1, \ldots, n \\ d_{i,j} > L_{i,j}, i, j = 1, \ldots, n, j \neq i \end{array} \right\} \tag{3}$$

where $w = (x_1, \ldots, x_n, y_1, \ldots, y_n, \theta_1, \ldots, \theta_n, v_1, \ldots, v_n)' \in \mathbb{R}^{4n}$ is the state vector of longitudinal and lateral positions, orientations, and speeds of all $n$ vehicles.

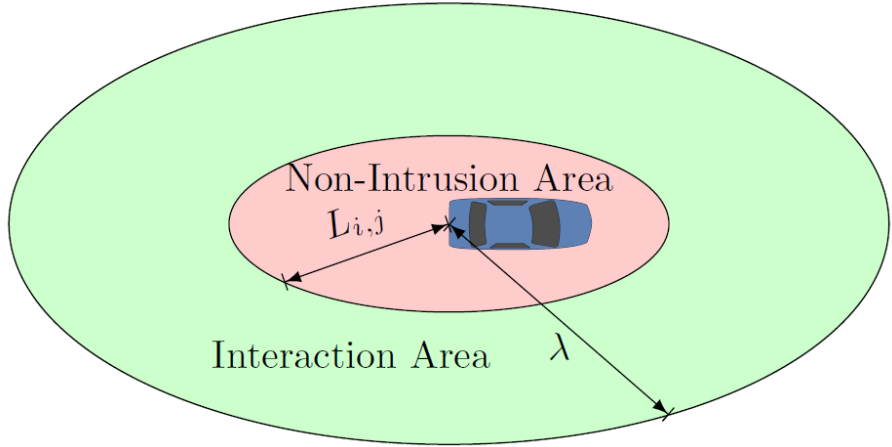

**Fig. 3:** Interaction and non-intrusion areas around a vehicle

The avoidance of collisions and road boundary violation can be achieved by using artificial potential functions. Let $V_{i,j}:(L_{i,j},+\infty)\to\mathbb{R}_+$, $U_i:(-a,a)\to\mathbb{R}_+$, $i,j=1,...,n$, $j\neq i$ be $C^2$ functions that satisfy the following properties

$$\lim_{d\to L_{i,j}^+}\left(V_{i,j}(d)\right)=+\infty \tag{4}$$

$$V_{i,j}(d)=0,\text{ for all }d\geq\lambda \tag{5}$$

$$V_{i,j}(d)=V_{j,i}(d),\ i,j=1,...,n,j\neq i \tag{6}$$

$$\lim_{y\to(-a)^+}\left(U_i(y)\right)=+\infty,\ \lim_{y\to a^-}\left(U_i(y)\right)=+\infty \tag{7}$$

$$U_i(0)=0 \tag{8}$$

where $\lambda$ is a positive constant that delineates the interaction area (Fig. 3), within which adjacent vehicles can in real-time exchange information (via sensor measurements or communication) such as relative positions and speeds. The families of functions $V_{i,j}$ and $U_i$ in (4), (5), (6), (7), and (8), are potential functions, which have been widely used to avoid collisions in mobile robotics (see [59], [83]). Note that the extent of the interaction area around the vehicle implies that, beyond repulsion by front vehicles (obstacles), repulsive forces are also exerted to the vehicle due to obstacles at its sides or behind it. The latter forces may also be called nudging forces to highlight that they do not exist in conventional anisotropic car-following by humans or even CAVs but are materializations of both TrafficFluid principles. No information is required for vehicles positioned out of the interaction area, as, in contrast with particle models from Mathematical Physics [96], the potential $V_{i,j}$ vanishes identically (see (5)) when the distance between two vehicles is greater than the interaction distance $\lambda$. Regarding road boundary respect, the employed artificial potential function $U_i$ tends to infinity when the vehicle reaches either boundary; and has virtually zero value at a distance from the boundaries (see (8)).

In view of all the above specifications and constraints, the CC design problem is a demanding problem, and several of the widely applied methodologies of Automatic Control [88], are not applicable for diverse reasons. Although other control design approaches may be possible, we address the problem via the use of Lyapunov functions, [61], and reach to a number of alternative CCs. This approach results in individual (decentralized) controllers for each vehicle with limited real-time information requirements about other neighboring vehicles. The derived nonlinear feedback CCs guarantee rigorously the fulfillment of all presented requirements and allow for a variety of different (parametric or functional) choices that entail accordingly different vehicle behaviors which all enjoy the same guaranteed qualitative properties.

Specifically, the family of CCs for vehicles operating on lane-free roads of constant width are constructed by means of Lyapunov functions expressed in terms of the mechanical energy of the vehicles and additional penalty terms that guarantee that all vehicle orientations remain within certain bounds. Two main approaches can be considered for the construction of the Lyapunov function. The first one is based on Newtonian mechanics, and the Lyapunov function is given by the formula

$$H(w) := \frac{1}{2}\sum_{i=1}^{n}\left(v_i \cos(\theta_i) - v^*\right)^2 + \frac{b}{2}\sum_{i=1}^{n} v_i^2 \sin^2(\theta_i) \\ + \sum_{i=1}^{n} U_i(y_i) + \frac{1}{2}\sum_{i=1}^{n}\sum_{j \neq i} V_{i,j}(d_{i,j}) + A\sum_{i=1}^{n}\left(\frac{1}{\cos(\theta_i) - \cos(\varphi)} - \frac{1}{1 - \cos(\varphi)}\right) \tag{9}$$

where $v^* \in (0, v_{\max})$ is the desired longitudinal speed of all vehicles, $b, A > 0$ are constants and

$$\cos(\varphi) > v^* / v_{\max}. \tag{10}$$

The first two terms ($\frac{1}{2}\sum_{i=1}^{n}\left(v_i \cos(\theta_i) - v^*\right)^2 + \frac{b}{2}\sum_{i=1}^{n} v_i^2 \sin^2(\theta_i)$) are related to the kinetic energy of the system of $n$ vehicles relative to an observer moving along the $x-$direction with speed equal to $v^*$ (as in classical mechanics); they penalize the deviation of the longitudinal and lateral speeds from their desired values $v^*$ and zero, respectively. The third and fourth terms ($\sum_{i=1}^{n} U_i(y_i)$ and $\frac{1}{2}\sum_{i=1}^{n}\sum_{j \neq i} V_{i,j}(d_{i,j})$), are based on the potential functions involved in (4) and (7) and are related to the potential energy of the system. While the Lyapunov function (9) has characteristics of a size function (see [89]), it is not a (global) size function, since $H$ takes finite values even when $v_i \notin (0, v_{\max})$. The following lemma shows the properties of the Lyapunov function $H$.

**Lemma 1:** *Let constants* $A > 0$, $v_{\max} > 0$, $v^* \in (0, v_{\max})$, $L_{i,j} > 0$, $i, j = 1,...,n$, $i \neq j$, $\lambda > 0$, $\varphi \in \left(0, \frac{\pi}{2}\right]$ *that satisfies (10), and define the function* $H : \Omega \to \mathbb{R}_+$ *by means of (9), where* $\Omega$ *is given by (3). Then, there exist non-decreasing functions* $\tau : \mathbb{R}_+ \to [0, \varphi)$, $\phi_i : \mathbb{R}_+ \to [0, a)$, $i = 1,...,n$, *and for each pair* $i, j \in \{1,...,n\}$, $i \neq j$, *there exist non-increasing functions* $\rho_{i,j} : \mathbb{R}_+ \to (L_{i,j}, \lambda]$ *with* $\rho_{i,j}(s) \equiv \rho_{j,i}(s)$, *such that the following implication holds:*

$$w \in \Omega \Rightarrow |\theta_i| \leq \tau(H(w)), |y_i| \leq \phi_i(H(w)), d_{i,j} \geq \rho_{i,j}(H(w)), \text{ for } i, j = 1,...,n, j \neq i.$$

In the second approach, the Lyapunov function is inspired by relativistic mechanics where the speed limits of relativistic mechanics (± the speed of light) are replaced by the two speed limits for vehicles (zero and the speed limit $v_{max}$ of the road, respectively; see [97]). In this case, the Lyapunov function is given by the formula

$$H_R(w) := \frac{v_{max}^2}{2} \sum_{i=1}^{n} \frac{\left(v_i \cos(\theta_i) - v^*\right)^2 + b v_i^2 \sin^2(\theta_i)}{(v_{max} - v_i) v_i} + \sum_{i=1}^{n} U_i(y_i)$$
$$+ \frac{1}{2} \sum_{i=1}^{n} \sum_{j \neq i}^{n} V_{i,j}(d_{i,j}) + A \sum_{i=1}^{n} \left( \frac{1}{\cos(\theta_i) - \cos(\varphi)} - \frac{1}{1 - \cos(\varphi)} \right) \quad (11)$$

The Lyapunov function $H_R$ is also a size function (see [89]) preventing the state $w$ to leave the state space $\Omega$.

**Lemma 2**: Let constants $A > 0$, $v_{max} > 0$, $v^* \in (0, v_{max})$, $L_{i,j} > 0$, $i, j = 1,...,n$, $i \neq j$, $\lambda > 0$, $\varphi \in \left(0, \frac{\pi}{2}\right]$ that satisfies (10), and define the function $H_R : \Omega \to \mathbb{R}_+$ by means of (11), where $\Omega$ is given by (3). Then, there exist non-decreasing functions $\phi_i : \mathbb{R}_+ \to [0, a)$, $i = 1,...,n$, $\ell_2 : \mathbb{R}_+ \to [v^*, v_{max})$, $\tau : \mathbb{R}_+ \to [0, \varphi)$, a non-increasing function $\ell_1 : \mathbb{R}_+ \to (0, v^*]$, and, for each pair $i, j = 1,...,n$, $i \neq j$, there exist non-increasing functions $\rho_{i,j} : \mathbb{R}_+ \to (L_{i,j}, \lambda]$ with $\rho_{i,j}(s) \equiv \rho_{j,i}(s)$, such that the following implications hold:

$$w \in \Omega \Rightarrow \ell_1(H_R(w)) \leq v_i \leq \ell_2(H_R(w)), |\theta_i| \leq \tau(H_R(w)), |y_i| \leq \phi_i(H_R(w)), d_{i,j} \geq \rho_{i,j}(H_R(w))$$
$$\text{for } i, j = 1,...,n, j \neq i.$$

### 2.2.1 Newtonian Cruise Controller (NCC)

The CC obtained from the Lyapunov function in (9) is given by

$$\delta_i = \arctan\left( \frac{\sigma_i}{v_i} \left( Z_i(w) - U_i'(y_i) - \sum_{j \neq i} p_{i,j} V_{i,j}'(d_{i,j}) \frac{(y_i - y_j)}{d_{i,j}} - b \sin(\theta_i) F_i \right) \right.$$
$$\left. \times \left( v^* + \frac{A}{v_i \left( \cos(\theta_i) - \cos(\varphi) \right)^2} + v_i \cos(\theta_i)(b-1) \right)^{-1} \right) \quad (12)$$

$$F_i = -\frac{1}{\cos(\theta_i)} \left( k_i(w) \left( v_i \cos(\theta_i) - v^* \right) + \Lambda_i(w) \right) \quad (13)$$

where

$$k_i(w) = \gamma + \frac{\Lambda_i(w)}{v^*} + \frac{v_{max} \cos(\theta_i)}{v^*(v_{max} \cos(\theta_i) - v^*)} \ell\left(-\Lambda_i(w)\right) \quad (14)$$

$$Z_i(w) := -\Gamma v_i \sin(\theta_i) + \sum_{j \neq i} \kappa_{i,j}(d_{i,j}) \left( g(v_j \sin(\theta_j)) - g(v_i \sin(\theta_i)) \right) \quad (15)$$

$$\Lambda_i(w) = \sum_{j \neq i} V'_{i,j}(d_{i,j}) \frac{(x_i - x_j)}{d_{i,j}} - \sum_{j \neq i} \kappa_{i,j}(d_{i,j}) \big( g(v_j \cos(\theta_j)) - g(v_i \cos(\theta_i)) \big) \tag{16}$$

$\gamma, \Gamma > 0$ are constants, $\ell \in C^1(\mathbb{R})$, $g \in C^1(\mathbb{R})$, and $\kappa_{i,j} \in C^1\big((L_{i,j}, +\infty), \mathbb{R}_+\big)$ for $i, j = 1, ..., n$, $j \neq i$ are functions that satisfy

$$\max(0, x) \leq \ell(x), \text{ for all } x \in \mathbb{R} \tag{17}$$
$$g'(x) > 0, \text{ for all } x \in \mathbb{R} \tag{18}$$
$$\kappa_{i,j}(d) = 0, \text{ for all } d \geq \lambda \tag{19}$$
$$\kappa_{i,j}(d) = \kappa_{j,i}(d), \ i, j = 1, ..., n, \ j \neq i \tag{20}$$

Using (1), (12), (13), it can be shown that the derivative of the Lyapunov function, satisfies $\nabla H(w) \dot{w} \leq 0$. Then, by exploiting Barbălat's Lemma ([20]), we further get the following theorem (see [54], [53]):

**Theorem 1**: *For every $w_0 \in \Omega$ there exists a unique solution $w(t) \in \Omega$ of the initial-value problem (1), (12), (13) with initial condition $w(0) = w_0$. The solution $w(t) \in \Omega$ is defined for all $t \geq 0$ and satisfies, for $i = 1, ..., n$,*

$$\lim_{t \to +\infty} (v_i(t)) = v^*, \ \lim_{t \to +\infty} (\theta_i(t)) = 0$$

$$\lim_{t \to +\infty} (u_i(t)) = 0, \ \lim_{t \to +\infty} (F_i(t)) = 0.$$

*Moreover, there exist non-decreasing functions $\bar{F} : \mathbb{R}_+ \to \mathbb{R}_+$ and $\bar{u} : \mathbb{R}_+ \to \mathbb{R}_+$ such that the following inequalities hold for every solution $w(t) \in \Omega$ of (1), (12), (13):*

$$|F_i(t)| \leq \bar{F}(H(w(0))), \ |\delta_i(t)| \leq \bar{u}(H(w(0))), \text{ for all } t \geq 0, \ i = 1, ..., n.$$

Theorem 1 guarantees the speeds of all vehicles tend asymptotically to the given speed set-point, i.e., $\lim_{t \to +\infty}(v_i(t)) = v^*$, $i = 1, ..., n$; the orientations of all vehicles tend asymptotically to be parallel to the road, i.e., $\lim_{t \to +\infty}(\theta_i(t)) = 0$, $i = 1, ..., n$; the accelerations and the steering angles tend asymptotically to zero, i.e., $\lim_{t \to +\infty}(F_i(t)) = \lim_{t \to +\infty}(\delta_i(t)) = 0$, $i = 1, ..., n$. While the bounds on acceleration and steering angle are determined by the technical characteristics of the vehicles and the road, the last two inequalities in Theorem 1 indicate that these control inputs are restricted by the function $H$ with finite energy (i.e., $H(w(0))$). This limitation is acceptable, since it is not possible to avoid collisions for any arbitrary initial condition under the presence of technical constraints.

The NCC has terms involving the functions $\kappa_{i,j}$ that imitate the action of viscosity in fluid flow, and can give faster convergence of the speeds. When the viscous terms $\kappa_{i,j}$ are present, the controllers require, in addition to relative position measurements, also real-time measurements of relative speeds of neighboring vehicles. Note that the viscosity does not act as a friction term but rather assists in equalizing the speeds of neighboring vehicles. The function $k_i(w)$ is a state-dependent controller gain, appropriately designed to ensure that the speed of each vehicle respects the speed limits. Notice that the properties in (5) and (19) guarantee that the NCC is decentralized (per vehicle) and requires measurement of the own state $(x_i, y_i, \theta_i, v_i)$ and the relative positions and relative speeds $(x_j - x_i, y_j - y_i, \dot{x}_i - \dot{x}_j, \dot{y}_i - \dot{y}_j)$ of all vehicles that are in distance $d_{i,j} = \sqrt{(x_i - x_j)^2 + p_{i,j}(y_i - y_j)^2}$ less than $\lambda$, i.e., within the interaction area.

### 2.2.2 Pseudo-Relativistic Cruise Controller (PRCC)

The PRCC that corresponds to the Lyapunov function in (11) is given by

$$\delta_i = \arctan\left(\frac{\sigma_i}{\beta(v_i, \theta_i)}\left(G_i(w) - U'_i(y_i) - \tilde{a}(v_i, \theta_i)F_i - \sum_{j \neq i} p_{i,j} V'_{i,j}(d_{i,j}) \frac{(y_i - y_j)}{d_{i,j}}\right)\right) \quad (21)$$

$$F_i = \frac{1}{q(v_i, \theta_i)}\left(R_i(w) - \sum_{j \neq i} V'_{i,j}(d_{i,j}) \frac{(x_i - x_j)}{d_{i,j}}\right) \quad (22)$$

where

$$G_i(w) = -\bar{f}\left(v_i \sin(\theta_i)\right) + \sum_{j \neq i} \kappa_{i,j}(d_{i,j})\left(g\left(v_j \sin(\theta_j)\right) - g\left(v_i \sin(\theta_i)\right)\right) \quad (23)$$

$$R_i(w) = -f\left(v_i \cos(\theta_i) - v^*\right) + \sum_{j \neq i} \kappa_{i,j}(d_{i,j})\left(g\left(v_j \cos(\theta_j)\right) - g\left(v_i \cos(\theta_i)\right)\right) \quad (24)$$

and

$$q(v, \theta) := v_{\max}^2 \frac{v_{\max} v \cos(\theta) + v^* v_{\max} - 2v^* v}{2(v_{\max} - v)^2 v^2} \quad (25)$$

$$\beta(v, \theta) := \frac{A}{(\cos(\theta) - \cos(\varphi))^2} + v_{\max}^2 \frac{(b-1)v\cos(\theta) + v^*}{(v_{\max} - v)} \quad (26)$$

$$\tilde{a}(v, \theta) := \frac{b v_{\max}^3 \sin(\theta)}{2(v_{\max} - v)^2 v} \quad (27)$$

The $C^1$ functions $f:(-v^*, v_{\max} - v^*) \to \mathbb{R}$, $\bar{f}:(-v_{\max}, v_{\max}) \to \mathbb{R}$ involved in (23) and (24) satisfy

$$x f(x) > 0, \; x \bar{f}(x) > 0, \text{ for } x \neq 0 \quad (28)$$

and $g \in C^1(\mathbb{R})$ satisfies (18).

Compared to the NCC, the PRCC is simpler, since it does not use state-dependent controller gains to restrict the speed in $(0, v_{\max})$. Also, notice that (5) and (19) guarantee again that the PRCC is decentralized (per vehicle) and require measurement of the own state $(x_i, y_i, \theta_i, v_i)$ and the relative positions and relative speeds $(x_j - x_i, y_j - y_i, \dot{x}_i - \dot{x}_j, \dot{y}_i - \dot{y}_j)$ of all vehicles that are within the interaction area.

Again, it should be noted that, when the PRCC is inviscid ($\kappa_{i,j}(d) \equiv 0$), then it does not require real-time measurements of the speeds of the adjacent vehicles. The PRCC has the same qualitative properties as the NCC (see [53]):

**Theorem 2:** *For every $w_0 \in \Omega$ there exists a unique solution $w(t) \in \Omega$ of the initial-value problem (1), (21), (22) with initial condition $w(0) = w_0$. The solution $w(t) \in \Omega$ is defined for all $t \geq 0$ and satisfies, for $i = 1, \ldots, n$,*

$$\lim_{t \to +\infty}(v_i(t)) = v^*, \ \lim_{t \to +\infty}(\theta_i(t)) = 0$$

$$\lim_{t \to +\infty}(u_i(t)) = 0, \ \lim_{t \to +\infty}(F_i(t)) = 0.$$

*Moreover, there exist non-decreasing functions $\tilde{F} : \mathbb{R}_+ \to \mathbb{R}_+$ and $\tilde{u} : \mathbb{R}_+ \to \mathbb{R}_+$ such that the following inequalities hold for every solution $w(t) \in \Omega$ of 1, 21, 22:*

$$|F_i(t)| \leq \tilde{F}(H(w(0))), \ |\delta_i(t)| \leq \tilde{u}(H(w(0))), \text{ for all } t \geq 0, \ i = 1, \ldots, n$$

### 2.2.3 Properties of the NCC and PRCC and lane-based operation

System (1) under the change of coordinates

$$\tilde{x}_i(t) = x_i(t) - \int_0^t \left( n^{-1} \sum_{j=1}^n v_j(s) \cos(\theta_j(s)) \right) ds - n^{-1} \sum_{j=1}^n x_j(0)$$

is given by (after a slight abuse of notation – dropping the tildes),

$$\begin{aligned}
\dot{x}_i &= v_i \cos(\theta_i) - n^{-1} \sum_{j=1}^n v_j \cos(\theta_j) \\
\dot{y}_i &= v_i \sin(\theta_i) \\
\dot{\theta}_i &= \frac{\sigma_i}{v_i} \tan(\delta_i) \\
\dot{v}_i &= F_i
\end{aligned} \quad (29)$$

The system in (29) evolves on the set $\Omega$ defined by (3) and satisfies $\frac{d}{dt}\left( \sum_{i=1}^n x_i \right) = 0$; hence it satisfies the additional condition $\sum_{i=1}^n x_i = 0$ (center of mass transformation; we have brought the $x$-coordinate of the center of mass of the vehicles to be at 0). The set of equilibrium points of the transformed closed-loop system of (29) with the NCC or PRCC is described by

$$E = \left\{ w \in \Omega : \begin{array}{c} v_i - v^* = \theta_i = \sum_{j \neq i} V'_{i,j}(d_{i,j}) \dfrac{(x_i - x_j)}{d_{i,j}} = 0, i = 1,...,n \\ U'_i(y_i) + \sum_{j \neq i} p_{i,j} V'_{i,j}(d_{i,j}) \dfrac{(y_i - y_j)}{d_{i,j}} = 0, i = 1,...,n \\ \sum_{i=1}^{n} x_i = 0 \end{array} \right\}.$$

For every solution $w(t)$ of the closed-loop system in (29) with the NCC or PRCC, it holds that $\lim_{t \to +\infty} (\dot{w}(t)) = 0$. If for some initial condition $w_0 \in \Omega$ with $\sum_{i=1}^{n} x_i(0) = 0$, the corresponding solution $w(t)$ is bounded, then there exists a compact set $\bar{K} \subset \Omega$ such that $w(t) \in \bar{K}$ for all $t \geq 0$. Consequently, in this case the omega limit set $\omega(w_0)$ is non-empty and satisfies $\omega(w_0) \subseteq E$. Since $\lim_{t \to +\infty} (dist(w(t), \omega(w_0))) = 0$ and $\omega(w_0) \subseteq E$ (which implies that $dist(w(t), \omega(w_0)) \geq dist(w(t), E)$) we can guarantee in this case that $\lim_{t \to +\infty} (dist(w(t), E)) = 0$. Thus, the NCC and PRCC guarantee that every bounded solution of the closed-loop system tends to the invariant set $E$ as $t \to +\infty$. However, it should be noted that the limit $\lim_{t \to +\infty} (dist(w(t), E)) = 0$ does not imply the existence of a point $w^* \in E$ for which $\lim_{t \to +\infty} (|w(t) - w^*|) = 0$. The latter implies that we cannot guarantee that an "ultimate" arrangement exists for the vehicles. On the other hand, simulations have shown that in all cases that we have tested, the solutions of the closed-loop system in (29) with the NCC or the PRCC are bounded and there exists a point $w^* \in E$ (that depends on the initial condition) for which $\lim_{t \to +\infty} (|w(t) - w^*|) = 0$, i.e., the vehicles tend to a final arrangement.

The design of the NCC and the PRCC and the proof of the convergence of the speeds to the desired speed $v^*$ are based on Barbălat's lemma, and thus, the convergence is only asymptotic and not exponential. On the other hand, when we consider the NCC for lane-based driving, it is possible to show asymptotic convergence of the solutions to the set of equilibrium points of the system; and, under certain conditions, we can also show exponential convergence. Analogous results also hold for the PRCC. The design of CCs for the case of a ring-road was also reported in [94].

The resulting closed-loop system for $n$ identical vehicles on a straight line (lane-based case) under the inviscid NCC with same minimum safety distance $L$ and $\lambda < 2L$ is given by the following set of ODEs

$$\begin{aligned} \dot{s}_i &= v_{i-1} - v_i \quad , \quad i = 2,...,n \\ \dot{v}_1 &= -k_1(s_2)(v_1 - v^*) - V'(s_2) \\ \dot{v}_i &= -k_i(s_i, s_{i+1})(v_i - v^*) + V'(s_i) - V'(s_{i+1}) \quad , \quad i = 2,...,n-1 \\ \dot{v}_n &= -k_n(s_n)(v_n - v^*) + V'(s_n) \end{aligned} \tag{30}$$

where

$$\begin{aligned}
s_i &= x_{i-1} - x_i \quad , \quad i = 2,\ldots,n \\
k_1(s_2) &= \gamma + \ell\left(-V'(s_2)\right) \\
k_i(s_i, s_{i+1}) &= \gamma + \ell\left(V'(s_i) - V'(s_{i+1})\right) \quad , \quad i = 2,\ldots,n-1 \\
k_n(s_n) &= \gamma + \ell\left(V'(s_n)\right)
\end{aligned} \tag{31}$$

$\gamma > 0$, $v^* \in (0, v_{\max})$ are constants, $\ell$ satisfies (17) and $V \in C^2((L,+\infty); \mathbb{R}_+)$ is a potential function that satisfies

$$\begin{aligned}
&\lim_{s \to L^+}(V(s)) = +\infty, \\
&V''(s) \geq 0, \\
&V(s) = 0, \quad \text{for } s \geq \lambda \\
&V'(s) < 0, \quad \text{for } L < s < \lambda
\end{aligned} \tag{32}$$

Notice that the state space of the model (30) is given by

$$\tilde{\Omega} = \left\{(s_2,\ldots,s_n, v_1,\ldots,v_n) \in \mathbb{R}^{2n-1} : \min_{i=2,\ldots,n}(s_i) > L, \max_{i=1,\ldots,n}(v_i) \leq v_{\max}, \min_{i=1,\ldots,n}(v_i) \geq 0\right\} \tag{33}$$

and that the set of equilibrium points of (30) is given by

$$\tilde{E} = \left\{(s_2,\ldots,s_n, v_1,\ldots,v_n) \in \mathbb{R}^{2n-1} : \min_{i=2,\ldots,n}(s_i) \geq \lambda, v_i = v^*, i = 1,\ldots,n\right\} \subset \tilde{\Omega} \tag{34}$$

where $v^* \in (0, v_{\max})$ is the desired speed.

Recall, that for the NCC in the lane-free case, we could only show asymptotic convergence of the speeds to the desired speed by means of Barbălat's Lemma. For the lane-based case however, using the Lyapunov function

$$\tilde{H}(s, v) = \frac{1}{2}\sum_{i=1}^{n}(v_i - v^*)^2 + \sum_{i=2}^{n}V(s_i)$$

and LaSalle's Invariance principle, [34], the following theorem holds (see [57])

**Theorem 3:** *For every initial condition $(s(0), v(0)) \in \tilde{\Omega}$, the solution $(s(t), v(t)) \in \tilde{\Omega}$ of (30) is defined for all $t \geq 0$ and satisfies*

$$s_i(t) \leq \max(\lambda, s_i(0)) + \gamma^{-1} v_{\max}, \text{ for all } t \geq 0 \text{ and } i = 2,\ldots,n. \tag{35}$$

*Moreover, $\lim_{t \to +\infty}(v_i(t)) = v^*$ for all $i = 1,\ldots,n$, $\lim_{t \to +\infty}(V(s_i(t))) = 0$, $\liminf_{t \to +\infty}(s_i(t)) \geq \lambda$ for all $i = 2,\ldots,n$, and $\lim_{t \to +\infty}\left(\text{dist}\left((s(t), v(t)), \tilde{E}\right)\right) = 0$. for all $i = 2,\ldots,n$.*

When the vehicles have large initial distances between them, then they simply adjust their speeds without affecting each other.

**Proposition 1:** *Suppose that* $s_i(0) \geq \max\left(\lambda - \bar{\omega}^{-1}(v_{i-1}(0) - v_i(0)), \lambda\right)$ *for* $i = 2,...,n$, *where* $\bar{\omega} = \gamma + \ell(0)$. *Then the solution of the model (30) is given by the equations:*

$$\begin{aligned}
v_i(t) &= v^* + \exp(-\bar{\omega}t)(v_i(0) - v^*) \quad, \quad i = 1,...,n \\
s_i(t) &= s_i(0) + \bar{\omega}^{-1}(v_{i-1}(0) - v_i(0))(1 - \exp(-\bar{\omega}t)) \quad, \quad i = 2,...,n
\end{aligned} \quad (36)$$

It should be noted again that LaSalle's principle does not guarantee uniform attraction to the set $\tilde{E}$, [34]. In order to be able to show uniform global attractivity properties for the set $\tilde{E}$, the following theorem provides a strict Lyapunov function for system (30).

**Theorem 4:** *For every* $\bar{\beta} > 0$, *there exist non-decreasing functions* $R \in C^1(\mathbb{R}_+;(0,+\infty))$, $\kappa \in C^0(\mathbb{R}_+;(0,+\infty))$ *such that the following inequalities hold for all* $(s,v) \in \tilde{\Omega}$:

$$H(s,v) \leq W(s,v) \leq \kappa(H(s,v))H(s,v)$$

$$\dot{W}(s,v) \leq -\gamma\bar{\beta}\sum_{i=1}^{n}(v_i - v^*)^2 - \frac{1}{8}\sum_{i=2}^{n}4^i(V'(s_i))^2$$

*where* $W : \tilde{\Omega} \to \mathbb{R}_+$ *is defined by the equation*

$$W(s,v) := R(H(s,v))H(s,v) - \sum_{i=2}^{n}4^i V'(s_i)(v_i - v^*) \text{, for all } (s,v) \in \tilde{\Omega} \quad (37)$$

*and* $\dot{W}(s,v)$ *denotes the time derivative of* $W$ *along the solutions of (30).*

In [57], using the Lyapunov function (37), we establish uniform global asymptotic stability with respect to the measures $dist((s,v),\tilde{E})$ and $W(s,v)$ (see [93] for the notion of asymptotic stability with respect to two measures).

**Theorem 5:** *There exist a function* $\psi \in KL$ *and a function* $\tilde{z} \in K_\infty$ *such that for every initial condition* $(s(0),v(0)) \in \tilde{\Omega}$ *the solution* $(s(t),v(t)) \in \tilde{\Omega}$ *of (30) is defined for all* $t \geq 0$ *and satisfies*

$$\tilde{z}\left(dist((s(t),v(t)),\tilde{E})\right) \leq W(s(t),v(t)) \leq \psi\left(W(s(0),v(0)),t\right) \text{, for all } t \geq 0.$$

## 2.3 Generalized Driving Situations

We next provide significant generalizations of the CC design problem, which open the way to consideration of many real complex driving situations, including the possible presence of on-ramps and off-ramps. Specifically, we allow for:
- Consideration of roads with variable width, which is important, among others, for implementation of internal boundary control, an innovative and highly efficient traffic control measure facilitated by lane-free driving, [70].

- Each vehicle may move within its own (curved) corridor boundaries, that may or may not coincide with the actual road boundaries. This allows to tackle many practical situations, including swift vehicle merging from on-ramps onto the mainstream of a highway or vehicles exiting at off-ramps.
- Possibility to have different individual desired speeds for each vehicle.

We use the bicycle kinematic model and we suppose that each vehicle $i = 1,...,n$ follows a path which is contained in the set (corridor)

$$Y_i = \{(x, y) : x \in \mathbb{R}, \alpha_i(x) < y < \beta_i(x)\} \tag{38}$$

where $\alpha_i : \mathbb{R} \to \mathbb{R}$ and $\beta_i : \mathbb{R} \to \mathbb{R}$ are $C^2$ functions with the following property: there exist constants $\varphi \in \left[0, \frac{\pi}{2}\right)$ (the maximum angle of orientation of each vehicle), $r_{max} > r_{min} > 0$ and $\gamma_{max}, \gamma_{min} \in \mathbb{R}$ with $\gamma_{max} > \gamma_{min}$ such that

$$\gamma_{max} \geq \beta_i(x), \ \alpha_i(x) \geq \gamma_{min} \ \text{for all} \ x \in \mathbb{R}, \ i = 1,...,n \tag{39}$$

$$r_{max} \geq \beta_i(x) - \alpha_i(x) \geq r_{min} > 0 \ \text{for all} \ x \in \mathbb{R}, \ i = 1,...,n \tag{40}$$

$$\max\left(\sup_{x \in \mathbb{R}}(|\alpha_i'(x)|), \sup_{x \in \mathbb{R}}(|\beta_i'(x)|)\right) < \tan(\varphi) \ \text{for all} \ x \in \mathbb{R}, \ i = 1,...,n \tag{41}$$

$$\sup_{x \in \mathbb{R}}(|\alpha_i''(x)| + |\beta_i''(x)|) < +\infty \ \text{for all} \ x \in \mathbb{R}, \ i = 1,...,n. \tag{42}$$

The inequality in (41) determines the maximum rate of change of each corridor boundary and ensures that the rate of change will not be higher than the maximum angle of orientation of each vehicle. The inequalities (40) determine the minimum and maximum width of each corridor $Y_i$, $i = 1,...,n$. Note that the corridor $Y_i$ of a vehicle may or may not have common points with the corridor $Y_j$ (with $j \neq i$) of a different vehicle. When all vehicles share the same corridor $Y_i$, namely, when $\alpha_i(x) = \alpha(x)$ and $\beta_i(x) = \beta(x)$, for all $i = 1,...,n$, then, the curves $\alpha(x)$ and $\beta(x)$ simply describe the right and left boundary of the road, respectively. Let $v_i^* \in (0, v_{max})$ be the desired speed of each vehicle $i = 1,...,n$ and assume that

$$\cos(\varphi) > \max\left\{\frac{\max_{i=1,...,n}\{v_i^*\}}{v_{max}}, \frac{1}{3}\right\} . \tag{43}$$

Let $\sigma : \mathbb{R} \to (0,1]$ be a non-increasing $C^2$ function that satisfies

$$\sigma(x) = 1 \ \text{for} \ x \leq \varepsilon \ \text{and} \ \sigma(x)x \leq M \ \text{for} \ x \geq \varepsilon \tag{44}$$

$$\sup_{x \in \mathbb{R}}(|\sigma'(x)| + |\sigma''(x)|) < +\infty \tag{45}$$

where $M$ and $\varepsilon$ are positive constants. Let $\bar{b} > 0$ be a constant and let $V_{i,j} : (L_{i,j}, +\infty) \to \mathbb{R}_+$, $U_i : (-1,1) \to \mathbb{R}_+$, $i, j = 1,...,n$, $j \neq i$ be $C^2$ functions that satisfy the properties shown in (4) – (8) with $a = 1$. Define

$$\bar{H}(w) = \frac{v_{\max}^2}{2} \sum_{i=1}^{n} \frac{(v_i \cos(\theta_i) - \zeta_i(w))^2}{(v_{\max} - v_i)v_i} + \frac{\bar{b}v_{\max}^2}{2} \sum_{i=1}^{n} \frac{(\sin(\theta_i) - g_i(x_i, y_i)\cos(\theta_i))^2}{\cos(\theta_i) - \cos(\varphi)}$$
$$+ \sum_{i=1}^{n} U_i \left( \frac{2y_i - (\beta_i(x_i) + \alpha_i(x_i))}{\beta_i(x_i) - \alpha_i(x_i)} \right) + \frac{1}{2} \sum_{i=1}^{n} \sum_{j \neq i} V_{i,j}(d_{i,j}) \tag{46}$$

where

$$g_i(x, y) := \frac{(\beta_i(x) - y)\alpha_i'(x) + (y - \alpha_i(x))\beta_i'(x)}{\beta_i(x) - \alpha_i(x)}, \text{ for all } (x, y) \in Y_i \tag{47}$$

and $\zeta_i$, $i = 1, \ldots, n$, are $C^1$ functions defined by

$$\zeta_i(w) = v_i^* \sigma\left(\Phi_i(w) + g_i(x_i, y_i)\Xi_i(w)\right), \ w \in \bar{\Omega}, \ i = 1, \ldots, n \tag{48}$$

with $\Phi_i(w)$ and $\Xi_i(w)$, for $i = 1, \ldots, n$, given by

$$\Phi_i(w) := \sum_{j \neq i} V_{i,j}'(d_{i,j}) \frac{(x_i - x_j)}{d_{i,j}}, \text{ for } w \in \bar{\Omega} \tag{49}$$

$$\Xi_i(w) = \sum_{j \neq i} p_{i,j} V_{i,j}'(d_{i,j}) \frac{(y_i - y_j)}{d_{i,j}}, \text{ for } w \in \bar{\Omega} \tag{50}$$

and

$$\bar{\Omega} := \left\{ w \in \mathbb{R}^{4n} : \begin{array}{l} (x_i, y_i) \in Y_i, i = 1, \ldots, n \\ v_i \in (0, v_{\max}), |\theta_i| < \varphi, i = 1, \ldots, n \\ d_{i,j} > L_{i,j}, i, j = 1, \ldots, n, j \neq i \end{array} \right\}$$

The design of CCs for this generalized case is based on the Lyapunov-like function $\bar{H}$ defined by (46). This function has differences from the Lyapunov functions of (11), since here we consider the general case where each vehicle has its own desired speed $v_i^*$ and its own corridor $Y_i$. Notice that the Lyapunov-like function follows a pseudo-relativistic setting where the kinetic energy tends to infinity when the speed tends to zero or the speed limit $v_{\max}$.

The CCs that are designed based on the Lyapunov-like function $\bar{H}$ are termed Generalized Cruise Controllers (GCCs).

### 2.3.1 Generalized Cruise Controllers (GCCs)

Using the Lyapunov-like function in (47) we obtain the following GCC

$$F_i = -\frac{1}{\bar{q}(v_i, \theta_i, \zeta_i(w))} \left( \gamma \left(v_i \cos(\theta_i) - \zeta_i(w)\right) + \Phi_i(w) + g_i(x_i, y_i)\Xi_i(w) \right)$$
$$+ \frac{v_{\max}^2 \left( \tilde{Z}_i(w) + \sigma_i \sin(\theta_i) \tan(\delta_i) \right)}{\bar{q}(v_i, \theta_i, \zeta_i(w)) v_i (v_{\max} - v_i)} \tag{51}$$

and

$$\delta_i = \arctan\left(\frac{2\sigma_i\left(\cos(\theta_i) - \cos(\varphi)\right)^2}{\bar{b}v_{\max}^2 v_i h_i(x_i, y_i, \theta_i)}\left(-\Gamma v_i^2\left(\sin(\theta_i) - g_i(x_i, y_i)\cos(\theta_i)\right)\right.\right.$$
$$\left.\left. + \bar{b}v_{\max}^2 \frac{\sigma_i \cos(\theta_i)\xi_i(x_i, y_i, \theta_i)}{\cos(\theta_i) - \cos(\varphi)} - U_i'\left(\frac{2y_i - (\beta_i(x_i) + \alpha_i(x_i))}{\beta_i(x_i) - \alpha_i(x_i)}\right)\frac{2\sigma_i}{\beta_i(x_i) - \alpha_i(x_i)} - \sigma_i\Xi_i(w)\right)\right) \quad (52)$$

where $\gamma, \Gamma > 0$ are constants and

$$\tilde{Z}_i(w) := v_i^* \sigma'\left(\Phi_i(w) + g_i(x_i, y_i)\Xi_i(w)\right)\frac{d}{dt}\left(\Phi_i(w) + g_i(x_i, y_i)\Xi_i(w)\right), \ i = 1, \ldots, n, \ w \in \bar{\Omega} \quad (53)$$

$$\bar{q}(v, \theta, \zeta) = v_{\max}^2 \frac{v_{\max} v \cos(\theta) + \zeta v_{\max} - 2\zeta v}{2(v_{\max} - v)^2 v^2}, \text{ for all } v \in (0, v_{\max}), \theta \in (-\varphi, \varphi), \zeta \in \mathbb{R} \quad (54)$$

$$h_i(x, y, \theta) = 2\cos(\theta)\left(\cos(\theta) - \cos(\varphi)\right) + \sin^2(\theta) + \sin(\theta)g_i(x, y)(\cos(\theta) - 2\cos(\varphi)),$$
$$\text{for all } (x, y) \in Y_i, \ \theta \in (-\varphi, \varphi), \ i = 1, \ldots, n \quad (55)$$

$$\xi_i(x_i, y_i, \theta) = \cos(\theta)\left(\frac{y_i - \alpha_i(x_i)}{\beta_i(x_i) - \alpha_i(x_i)}\beta_i''(x_i) + \left(1 - \frac{y_i - \alpha_i(x_i)}{\beta_i(x_i) - \alpha_i(x_i)}\right)\alpha_i''(x_i)\right)$$
$$+ \left((\beta_i(x_i) - \alpha_i(x_i))\sin(\theta) - g_i(x_i, y_i)\cos(\theta)\right)\frac{\beta_i'(x_i) - \alpha_i'(x_i)}{(\beta_i(x_i) - \alpha_i(x_i))^2}$$
$$\text{for all } (x_i, y_i) \in Y_i, \ \theta \in (-\varphi, \varphi), \ i = 1, \ldots, n. \quad (56)$$

The GCC for each vehicle requires measurements of its own state $(x_i, y_i, \theta_i, v_i)$ and the relative positions and relative speeds $(x_j - x_i, y_j - y_i, \dot{x}_i - \dot{x}_j, \dot{y}_i - \dot{y}_j)$ of all vehicles that are in distance $d_{i,j} = \sqrt{(x_i - x_j)^2 + p_{i,j}(y_i - y_j)^2}$ less than $\lambda$. Moreover, the GCC for each vehicle requires knowledge of the values of its own boundary functions $\alpha_i(x)$, $\beta_i(x)$ and their derivatives $\alpha_i'(x), \beta_i'(x)$ and $\alpha_i''(x), \beta_i''(x)$ at the current $x$-coordinate of the position of the vehicle. Therefore, the GCC is also decentralized. For the GCC, the following theorem was proved in [56].

**Theorem 6:** *For every $w_0 \in \bar{\Omega}$, the initial value problem (1), (51), (52) with initial condition $w(0) = w_0$ has a unique solution $w(t)$, defined for all $t \geq 0$, that satisfies $w(t) \in \bar{\Omega}$ for all $t \geq 0$.*

Theorem 6 guarantees that: (i) no collisions occur, (ii) no vehicle exits its movement corridor $Y_i = \{(x, y) : x \in \mathbb{R}, \alpha_i(x) < y < \beta_i(x)\}$, (iii) no vehicle has a speed greater than the speed limit, (iv) no vehicle moves in the upstream direction, and (v) no vehicle has orientation greater than $\varphi$. It should be noted that the GCC does not guarantee a positive lower bound for the speed of all vehicles, i.e., it may happen that $\lim_{t \to +\infty}(v_i(t)) = 0$ for some $i \in \{1, \ldots, n\}$, and the corresponding vehicle in this case would tend to stop.

When all vehicles are in constant-width corridors $Y_i$ defined by $\alpha_i(x) \equiv \alpha_i$, $\beta_i(x) \equiv \beta_i$ for $i = 1, \ldots, n$ with $\beta_i > \alpha_i$ (recall (40)), and all vehicles have the same speed set-points $v_i^* = v^*$ for $i = 1,...,n$, then, the GCC has the same qualitative properties as the NCC and PRCC. Namely, in this case the solution is not only well-defined, but in addition the speed $v_i$ of each vehicle converges (asymptotically) to the speed set-point $v^*$ and the orientation $\theta_i$, and accelerations of each vehicle converge (asymptotically) to 0 (see [56]).

**Theorem 7:** *Suppose that $\beta_i(x) \equiv \beta_i \in \mathbb{R}$ and $\alpha_i(x) \equiv \alpha_i \in \mathbb{R}$ with $r_{\max} \geq \beta_i - \alpha_i \geq r_{\min} > 0$ and $v_i^* = v^*$ for all $i = 1,...,n$. Then, for every solution $w(t)$ of the closed-loop system (1), (51), (52), the following equations hold for $i = 1,...,n$:*

$$\lim_{t \to +\infty} (v_i(t)) = v^*, \; \lim_{t \to +\infty} (\theta_i(t)) = 0$$

$$\lim_{t \to +\infty} (\Phi_i(w(t))) = \lim_{t \to +\infty} \left( U_i' \left( \frac{2y_i(t) - (\beta_i + \alpha_i)}{\beta_i - \alpha_i} \right) \frac{2}{\beta_i - \alpha_i} + \Xi_i(w(t)) \right) = 0$$

$$\lim_{t \to +\infty} (F_i(t)) = 0, \; \lim_{t \to +\infty} (u_i(t)) = 0$$

## 3. Derivation of Macroscopic Models from Microscopic Models

### 3.1 Particle Methods

Research works have used various methods for the formal derivation of macroscopic traffic flow models from microscopic models (see [5], [24], [37], [38], [40]). In particular, particle methods were used recently for the rigorous proof of existence of entropy solutions of traffic flow problems (see [39], [41]). Moreover, the convergence of the Navier-Stokes (NS) solutions produced by particle methods to weak solutions of 1-D compressible NS equations was shown in [51]. More specifically, it was shown in [51] that particle methods guarantee differential inequalities and conservation laws that hold for classical solutions of the macroscopic model, but there is no guarantee that the inequalities or the conservation laws hold (in a discretized form) for an arbitrary numerical approximation. This is a well-known problem (see [32], [48], [92] for the analysis in a finite-dimensional setting), which is closely related to the issue of numerical stability for the applied scheme. The success of particle methods has been restricted so far to flows with one spatial dimension. However, in traffic flow theory most macroscopic models involve only one spatial dimension. Consequently, particle methods seem to be the proper "tool" for the derivation of macroscopic traffic flow models.

Particle methods are also known as the method of "Smoothed Particle Hydrodynamics" (see [8], [24], [40], [41], [80], [99] and references therein) and have been used extensively for many difficult flow problems, as well as for Partial Differential Equations (PDEs) other than the NS equations (see [12], [22]).

It can be said that particle methods provide the link between the macroscopic fluid model and the microscopic ODE model that describes the interaction between the fluid particles. This is important from a mathematical viewpoint (Hilbert's 6[th] problem, [42]), even when the forces acting on the particles are fictitious. While fictitious forces have no place in Mathematical Physics, they may appear in the context of the design and use of CCs for automated vehicles. CCs manipulate the accelerations of the vehicles via algorithmic creation of forces that may be considered as fictitious or artificial from a physical point of view. Thus, the vehicles are considered as (self-driven) particles of the "traffic fluid".

All macroscopic models presented below are formally derived based on particle models. Adopting a methodology that is inspired from [24], [40], [41], we provide the macroscopic models for which the particle (microscopic) approximation yields the ODE models of vehicles under the CCs that are described above. The derivation was first given in [57] for constant-width roads. In all cases, we consider only one spatial coordinate (denoted by $x$), i.e., we neglect the lateral movement of the vehicles. This is certainly an accurate approximation for the lane-based operation of the CCs on a single lane. Regarding lane-free traffic, we believe that the derived macroscopic models present good approximations for the emerging traffic flow, and ongoing work aims to demonstrate this. The obtained macroscopic models are similar to the NS equations for the 1-D flow of a viscous compressible fluid with viscosity depending on the density of the fluid, [49], [50].

For the application of the particle method, we consider a set of $n$ particles (vehicles) of total mass $m > 0$. Each vehicle has mass $m/n$ and vehicle $i \in \{1,...,n\}$ is placed at $x_i \in \mathbb{R}$ with $n(x_i - x_{i+1}) > L$ for $i = 1,...,n-1$, where $L \geq 0$ is a constant. We define the inter-vehicle distance by

$$s_i = x_{i-1} - x_i \quad , \quad i = 2,...,n. \tag{57}$$

## 3.2 Traffic Flow of Human-Driven Vehicles

The first and most influential traffic flow model for human drivers was independently proposed by Lighthill and Whitham [68] and Richards [84] (LWR model). The well-known LWR model is a first-order model governed by one scalar hyperbolic PDE which is simply the continuity equation (conservation of mass)

$$\rho_t + (\rho v)_x = 0 \tag{58}$$

where $\rho > 0$ is the vehicle density and $v$ is the vehicle mean speed, and both are combined with an algebraic equation that determines the vehicle speed as a function of the vehicle density

$$v = F(\rho) \tag{59}$$

where $F:[0,+\infty) \to [0,+\infty)$ is a non-increasing function. The LWR model has been studied extensively in the literature (see for instance [40], [50]-[102] and it was shown in [24], by using particle methods, that its entropy solutions are appropriate limits, as $n \to +\infty$, of solutions of the ODE system

$$\dot{x}_1 = F(0), \; \dot{x}_i = F\left(\frac{m}{ns_i}\right), \; i = 2,...,n. \tag{60}$$

First-order traffic models, while well-studied, present several limitations, as they do not capture significant aspects of real traffic flow dynamics. Therefore, researchers proposed second-order models for traffic flow with human drivers which improve the modeling accuracy and some of the deficiencies of first-order models at the expense of higher complexity (an additional PDE). The most popular second-order model for human-driven vehicles is the ARZ model (see [4], [5], [109]), which consists of the continuity equation (58) and the following speed equation

$$v_t + (v + \rho F'(\rho))v_x = -\bar{k}(v - F(\rho)) \tag{61}$$

where $\bar{k} > 0$ is a constant. The connection of the ARZ model with a microscopic Follow-the-Leader model was shown in [4]. By following a particle methodology, it can be formally shown that the microscopic model that corresponds to the macroscopic model of (58) and (61) is the following system of ODEs

$$\dot{x}_i = v_i, \; i = 1,...,n \tag{62}$$

$$\dot{v}_i = mF'\left(\frac{m}{ns_i}\right)\frac{v_i - v_{i-1}}{ns_i^2} - \bar{k}\left(v_i - F\left(\frac{m}{ns_i}\right)\right), \; i = 2,...,n \tag{63}$$

with $v_1 = F(0)$.

Both the LWR and ARZ models have an important specific feature: they are anisotropic, [15]. This becomes apparent in the respective microscopic equations (60) and (63), where the speed or acceleration of vehicle $i$ depend on the state of vehicle $i-1$ but not on the state of vehicle $i+1$. Another common feature of the LWR and ARZ models is that both models are hyperbolic with no viscosity (or diffusion) terms. It should also be mentioned that the LWR model is immersed in the ARZ model in the following sense: if the initial condition for the ARZ model satisfies (59) everywhere, then the solution of the ARZ model coincides with the solution of the LWR model (see [31], [66]). This becomes apparent by defining

$$s = v - F(\rho) \tag{64}$$

which is a Riemann coordinate of the ARZ model and satisfies the equation:

$$s_t + vs_x = -\bar{k}s. \tag{65}$$

Both the LWR and the ARZ models have been extensively used in the literature for control purposes (see [8], [10], [64], [102], [104], [105], [110]). The recent book [106] contains many extensions of the ARZ model to various cases (multiple lanes, multiple vehicle classes) and explicit controller or observer designs for the linearization of the ARZ model for finite domains with boundary inputs.

It is important to note that the microscopic model of (62) and (63) coincides with the closed-loop system that is obtained in [91] when the vehicles follow a platoon formation and are under the effect of a nonlinear Follow-the-Leader CC, with one essential difference: the CC proposed in [55] requires a density-dependent relaxation constant, i.e., $\bar{k} > 0$ is not a constant but a function of $s_i$ ($\bar{k} = \bar{k}(s_i)$). Thus, the nonlinear Follow-the-Leader CC proposed in [55] gives rise to a macroscopic model that consists of the same equations as the ARZ model, with the difference that the relaxation coefficient is density-dependent, i.e., $\bar{k} = \bar{k}(\rho)$ in (61).

### 3.3 Macroscopic Models for CAVs

We consider the movement of $n$ identical vehicles with total mass $m > 0$ on a constant-width road under the PRCC of (21) and (22) or under the NCC of (12) and (13) or under the GCC of (51) and (52) when:
(i) the vehicles are constrained to move on the line $y = 0$ (longitudinal motion; $\theta_i(t) \equiv 0$),

(ii) there exist constants $\lambda > L > 0$ with $\lambda < 2L$, such that $V_{i,j}(s) = \Phi(ns)$ for all $i, j = 1,...,n$ and $s > L/n$, where $\Phi:(L,+\infty) \to \mathbb{R}_+$ is a $C^2$ function that satisfies $\lim_{d \to L^+}(\Phi(d)) = +\infty$ and $\Phi(d) = 0$ for all $d \geq \lambda$,

(iii) $\kappa_{i,j}(s) = n^2 K(ns)$, for all $i, j = 1,...,n$ and $s > L/n$, where $K:(L,+\infty) \to \mathbb{R}_+$ is a $C^1$ function that satisfies $K(d) = 0$ for all $d \geq \lambda$, and

(iv) the number of vehicles $n$ is very large (tends to infinity).

In all cases, the state space is the open set
$$\Omega_n = \left\{ (x_1,...,x_n, v_1,...,v_n) \in \mathbb{R}^n \times (0, v_{\max})^n : n(x_i - x_{i+1}) > L, i = 1,...,n-1 \right\}$$

### 3.3.1 Macroscopic Model for PRCCs (macro-PRCC)

The microscopic model of (1) under the PRCC of (21) and (22) is given by the following ODEs:
$$\dot{x}_i = v_i \quad , \quad i = 1, 2, ..., n \tag{66}$$
$$\tilde{q}(v_1)\dot{v}_1 = -f(v_1 - v^*) - n\Phi'(ns_2) + n^2 K(ns_2)(g(v_2) - g(v_1)) \tag{67}$$
$$\tilde{q}(v_i)\dot{v}_i = -f(v_i - v^*) + n(\Phi'(ns_i) - \Phi'(ns_{i+1}))$$
$$+ n^2 \left( K(ns_i)(g(v_{i-1}) - g(v_i)) + K(ns_{i+1})(g(v_{i+1}) - g(v_i)) \right)$$
$$\text{for } i = 2,...,n-1 \tag{68}$$
$$\tilde{q}(v_n)\dot{v}_n = -f(v_n - v^*) + n\Phi'(ns_n) + n^2 K(ns_n)(g(v_{n-1}) - g(v_n)) \tag{69}$$
where $\tilde{q}:(0, v_{\max}) \to (0, +\infty)$ is given by
$$\tilde{q}(v) := v_{\max}^2 \frac{v_{\max} v - 2v^* v + v^* v_{\max}}{2(v_{\max} - v)^2 v^2}. \tag{70}$$

The macroscopic model that corresponds to the microscopic model described by (66), (67), (68), (69) is (58) and the following PDE
$$\tilde{q}(v)(v_t + vv_x) + \frac{P'(\rho)}{\rho}\rho_x = \frac{1}{\rho}(g'(v)\mu(\rho)v_x)_x - f(v - v^*) \tag{71}$$

that hold for $t > 0$ and $x \in I(t)$, where $I(t) \subseteq \mathbb{R}$ is an appropriate interval, with state constraints $0 < \rho(t,x) < \rho_{\max}$, $0 < v(t,x) < v_{\max}$ for all $t > 0$ and $x \in I(t)$, where
$$\rho_{\max} = \frac{m}{L} \tag{72}$$
$$\mu(\rho) = \frac{m^2}{\rho} K\left(\frac{m}{\rho}\right), \quad P(\rho) = -m\Phi'\left(\frac{m}{\rho}\right) \text{ for } \rho \in (0, \rho_{\max}). \tag{73}$$

The macroscopic model given by (58) and (71) will be called macro-PRCC. Notice that the assumptions for PRCC imply that
$$\int_{\bar{\rho}}^{\rho_{\max}} \rho^{-2} P(\rho) d\rho = +\infty, \tag{74}$$

$$\mu(\rho) = 0, \ P(\rho) = 0 \text{ for all } \rho \in (0, \bar{\rho}] \tag{75}$$

where

$$\bar{\rho} = \frac{m}{\lambda} \tag{76}$$

is the interaction density.

### 3.3.2 Macroscopic Model for NCCs (macro-NCC)

The microscopic model of (1) under the NCC of (12) and (13) is given by the ODEs in (66) and the following ODEs:

$$\dot{v}_1 = -J(\tilde{G}_1)(v_1 - v^*) + \tilde{G}_1 \tag{77}$$

$$\dot{v}_i = -J(\tilde{G}_i)(v_i - v^*) + \tilde{G}_i, \text{ for } i = 2,...,n-1 \tag{78}$$

$$\dot{v}_n = -J(\tilde{G}_n)(v_n - v^*) + \tilde{G}_n \tag{79}$$

where

$$J(s) = \gamma + \frac{v_{max}\ell(s)}{v^*(v_{max} - v^*)} - \frac{s}{v^*} \tag{80}$$

and

$$\begin{aligned}
\tilde{G}_1 &= -n\Phi'(ns_2) + n^2 K(ns_2)(g(v_2) - g(v_1)) \\
\tilde{G}_i &= n\Phi'(ns_i) - n\Phi'(ns_{i+1}) + n^2 K(ns_i)(g(v_{i-1}) - g(v_i)) \\
&\quad + n^2 K(ns_{i+1})(g(v_{i+1}) - g(v_i)), \quad i = 1,...,n \\
\tilde{G}_n &= n\Phi'(ns_n) + n^2 K(ns_n)(g(v_{n-1}) - g(v_n))
\end{aligned} \tag{81}$$

The macroscopic model that corresponds to the microscopic model described by (66), (77), (78) and (76) is given by the continuity equation (58) and the following speed equation for $t > 0$ and $x \in I(t)$, where $I(t) \subseteq \mathbb{R}$ is an appropriate interval

$$v_t + vv_x + \frac{P'(\rho)}{\rho}\rho_x = \frac{1}{\rho}\left(g'(v)\mu(\rho)v_x\right)_x - J(G)(v - v^*) \tag{82}$$

where

$$G = -\frac{P'(\rho)}{\rho}\rho_x + \frac{1}{\rho}\left(\mu(\rho)g'(v)v_x\right)_x \tag{83}$$

and $\rho_{max}, \mu(\rho), P(\rho)$ are given by (72) and (73). Again, the macroscopic model given by (58) and (82) will be called macro-NCC and is to be considered with state constraints $0 < \rho(t,x) < \rho_{max}$, $0 < v(t,x) < v_{max}$ for all $t > 0$ and $x \in I(t)$.

### 3.3.3 Macroscopic Model for GCCs (macro-GCC)

The microscopic model of (1) under the GCC of (51) and (52) is given by the ODEs in (66) and the following ODEs:

$$\frac{d}{dt}\left(\frac{v_1 - v^*\sigma\left(n\Phi'(ns_2)\right)}{c(v_1)}\right) = -\gamma \frac{c(v_1)}{v_{max}^2}\left(v_1 - v^*\sigma\left(n\Phi'(ns_2)\right) + n\gamma^{-1}\Phi'(ns_2)\right) \tag{84}$$

$$\frac{d}{dt}\left(\frac{v_i - v^*\sigma\left(n\left(\Phi'(ns_{i+1}) - \Phi'(ns_i)\right)\right)}{c(v_i)}\right) =$$

$$-\gamma \frac{c(v_i)}{v_{max}^2}\left(v_i - v^*\sigma\left(n\left(\Phi'(ns_{i+1}) - \Phi'(ns_i)\right)\right) + n\gamma^{-1}\left(\Phi'(ns_{i+1}) - \Phi'(ns_i)\right)\right)$$

$$\text{for } i = 2,...,n-1 \tag{85}$$

$$\frac{d}{dt}\left(\frac{v_n - v^*\sigma\left(-n\Phi'(ns_n)\right)}{c(v_n)}\right) = -\gamma \frac{c(v_n)}{v_{max}^2}\left(v_n - v^*\sigma\left(-n\Phi'(ns_n)\right) - n\gamma^{-1}\Phi'(ns_n)\right) \tag{86}$$

where

$$c(v) = \sqrt{(v_{max} - v)v}, \text{ for } v \in (0, v_{max}) \tag{87}$$

The macroscopic model that corresponds to the microscopic model described by (84), (85) and (86) is given by the continuity equation (58) and the following speed equation for $t > 0$ and $x \in I(t)$, where $I(t) \subseteq \mathbb{R}$ is an appropriate interval

$$Q(v,r)(v_t + vv_x) + \frac{P'(\rho)}{\rho}\rho_x = -\frac{v_{max}^2 v^*\sigma'(r)}{\rho(v_{max} - v)v}\left(\rho P'(\rho)v_x\right)_x - \gamma\left(v - v^*\sigma(r)\right) \tag{88}$$

and

$$r := \frac{P'(\rho)}{\rho}\rho_x \tag{89}$$

$$Q(v,r) := v_{max}^2 \frac{v_{max}v + (v_{max} - 2v)v^*\sigma(r)}{2(v_{max} - v)^2 v^2} \tag{90}$$

and $\rho_{max}, \mu(\rho), P(\rho)$ are given by (72) and (73). Again, the macroscopic model that consists of (58), (88) is named macro-GCC and is to be considered with state constraints $0 < \rho(t,x) < \rho_{max}$, $0 < v(t,x) < v_{max}$ for all $t > 0$ and $x \in I(t)$. Moreover, notice that if we define

$$z = \frac{v - v^*\sigma(r)}{\sqrt{(v_{max} - v)v}} \tag{91}$$

then we obtain the simplified equation

$$z_t + vz_x + \frac{c(v)}{v_{max}^2}\frac{P'(\rho)}{\rho}\rho_x = -\gamma \frac{c^2(v)}{v_{max}^2} z. \tag{92}$$

### 3.4 Relation Between Macroscopic and Microscopic Quantities

We next present the relations between the various parameters and functions involved in the microscopic models on one hand; and the corresponding macroscopic quantities involved in the macroscopic models on the other hand. Table 1 shows how all parameters and functions of the macroscopic models can be directly obtained from the corresponding microscopic models. Table 1 shows that by changing the functions and the parameters of the CCs, we can directly determine the physical properties of the ''traffic fluid''. In this sense, we may talk about an engineered or designed artificial fluid that approximates the actual emerging traffic flow. To appreciate the strength and importance of the correspondence between the CC settings and the characteristics of the emerging traffic fluid, it is instructive to notice that for general barotropic (or isentropic) flow of compressible Newtonian fluids (see [69]), the dynamic viscosity and the pressure are always increasing functions of the fluid density. However, for a traffic fluid that emerges from CAV driving, this may not be true (e.g., pressure or viscosity can have local minima) if we so decide. Thus, the traffic fluid can have very different physical properties from those of real compressible Newtonian fluids depending on the employed CCs.

It should also be mentioned here that for real gases there is no constant $\bar{\rho}$ for which (75) holds. Moreover, for general barotropic (or isentropic) flow of real gases, (74) holds with $\rho_{max} = +\infty$. On the other hand, for traffic fluids $\rho_{max}$ is finite and can be determined by the parameters of the CCs.

**Table 1: Relations between macroscopic and microscopic quantities**

|  | Macroscopic | Microscopic |
|---|---|---|
| Maximum Density | $\rho_{max}$ | $= \dfrac{m}{L}$ |
| Maximum Speed | $v_{max}$ | $= v_{max}$ |
| Desired Speed | $v^*$ | $= v^*$ |
| Interaction Density | $\bar{\rho}$ | $= \dfrac{m}{\lambda}$ |
| Dynamic Viscosity | $\mu(\rho)$ | $= \dfrac{m^2}{\rho} K\left(\dfrac{m}{\rho}\right)$ |
| Pressure | $P(\rho)$ | $= -m\Phi'\left(\dfrac{m}{\rho}\right)$ |

### 3.5 Discussion of Macroscopic Models for CAVs

The macroscopic models for CAVs, i.e., the macro-PRCC, for the macro-NCC, and the macro-GCC, are all fluid-like models. Many terms in the speed equations can be named as in classical fluid studies:

- The term $\dfrac{P'(\rho)}{\rho}\rho_x$ that appears in (71), (82) and (88) is the pressure term that also appears in NS equations.

- The term $\dfrac{1}{\rho}\left(g'(v)\mu(\rho)v_x\right)_x$ that appears in (71) and (82) is similar to the viscosity term $\dfrac{1}{\rho}\left(\mu(\rho)v_x\right)_x$ that also appears in the NS equations. The same holds when $P'(\rho) \geq 0$ for all $\rho \in (0, \rho_{\max})$ (i.e., when $P(\rho)$ is non-decreasing) for the term $-\dfrac{v_{\max}^2 v^* \sigma'(r)}{\rho(v_{\max}-v)v}(\rho P'(\rho) v_x)_x$ that appears in (88). In fact, this term is a viscosity term. However, note that in this case the viscosity term is induced by pressure and it is not an independent property of the system.

- The term $-f(v-v^*)$ that appears in (71), the term $-J(G)(v-v^*)$ that appears in (82) and the term $-\bar{\gamma}(v - v^* \sigma(r))$ that appears in (88), are all similar to friction terms or relaxation terms of the form $-\kappa(\rho,v)v$ that appear in the NS equations. However, it should be noted that the terms in the macroscopic models for CAVs describe the tendency of vehicles to adjust their speed to the given speed set-point $v^* \in (0, v_{\max})$, while the friction term in the NS equation has the tendency to bring the fluid speed to zero.

However, there are also some crucial differences between the macroscopic models and the 1-D compressible Navier-Stokes equations:
- In compressible fluid flow, the density is only required to be positive, while the speed can take arbitrary real values [69]. In contrast, in the macroscopic models that describe the evolution of traffic for vehicles under the proposed CCs, the density is required to be positive and less than a maximum density, while the speed is required to be positive and less than the speed limit of the road.
- In compressible fluid flow, the friction term tends to bring the fluid speed to zero; while in the obtained macroscopic models, the friction term tends to bring the speed to the given speed set-point $v^* \in (0, v_{\max})$.
- In compressible fluid flow, the pressure and viscosity terms are always present, and they never vanish; while in the derived macroscopic models, the pressure and viscosity terms vanish when the density is lower than a specific value (the interaction density).
- In the macro-PRCC and macro-GCC models, the acceleration is multiplied by a factor that blows up when the speed tends to zero or to the road speed limit. This factor is analogous to factors appearing in relativistic fluid mechanics and stems from the pseudo-relativistic character of the mechanical energy that was used for the design of PRCC and GCC.

The similarity of the macroscopic models for CAVs with the 1-D NS equations for a compressible fluid is shown in Table 2 (at the end of this manuscript), which reveals and highlights many aspects of the derived macroscopic models. The similarity of the macroscopic models for CAVs with the 1-D NS equations for a compressible fluid is not accidental. Contrary to the case of human-driven vehicles, the traffic induced by CAVs is isotropic and the vehicles behave as self-driven particles of a Newtonian fluid. The vehicles under the effect of the proposed CCs do not react to downstream vehicles only (as in conventional traffic), but also to upstream vehicles (the nudging effect).

Another important aspect of the derived macroscopic models for CAVs is the fact that they do not include non-local terms. There are some models in the literature for traffic flow of CAVs for which non-local terms play a crucial role (see [11], [13], [25], [30]) and can even capture possible nudging effects. For example, in [52] the non-local model that consists of (58) and the equation

$$v(t,x) = f_{NL}\left(\int_{x}^{x+\eta_{NL}} \omega_{NL}(s-x)\rho(t,s)ds\right) g_{NL}\left(\int_{x-\zeta_{NL}}^{x} \tilde{\omega}_{NL}(x-s)\rho(t,s)ds\right) \quad (93)$$

was proposed. In (93), the term $f_{NL}\left(\int_{x}^{x+\eta_{NL}} \omega_{NL}(s-x)\rho(t,s)ds\right)$ is the "look-ahead" term, where $\eta_{NL} > 0$ is a constant (reflecting the visibility area), $f_{NL}: \mathbb{R}_+ \to \mathbb{R}_+$ and $\omega_{NL}: \mathbb{R}_+ \to \mathbb{R}_+$ are non-increasing functions with $\int_{0}^{\eta_{NL}} \omega_{NL}(x)dx = 1$. Moreover, the term $g_{NL}\left(\int_{x-\zeta_{NL}}^{x} \tilde{\omega}_{NL}(x-s)\rho(t,s)ds\right)$ is the "look-behind" term that expresses the nudging effect, where $\zeta_{LN} > 0$ is a constant, $g_{LN}: \mathbb{R}_+ \to \mathbb{R}_+$ is a non-decreasing, bounded function, and $\tilde{\omega}_{LN}: \mathbb{R}_+ \to \mathbb{R}_+$ is a non-increasing function. It is clear that the non-local model that consists of (58) and (93) is very different from the fluid-like macroscopic models for CAVs, i.e., the macro-PRCC, the macro-NCC, and the macro-GCC.

The macro-NCC is a nonlinear system of PDEs that can be solved in some cases when the domain is $\mathbb{R}$ (i.e., when $x \in \mathbb{R}$). More specifically, we consider (58) and (82) with initial condition

$$(\rho[0], v[0]) = (\rho_0, v_0) \tag{94}$$

Define $\tilde{\omega} = \gamma + \dfrac{v_{max}\ell(0)}{v^*(v_{max} - v^*)}$. The following theorem guarantees that the macro-NCC has a classical solution that satisfies certain estimates.

**Theorem 8:** *Consider the initial-value problem (58), (82), (94), with $\rho_0 \in C^1(\mathbb{R}) \cap W^{1,\infty}(\mathbb{R})$, $v_0 \in C^2(\mathbb{R}) \cap W^{2,\infty}(\mathbb{R})$ that satisfies*

$$\inf_{x \in \mathbb{R}}(v_0'(x)) > -\tilde{\omega} \tag{95}$$

$$\sup_{x \in \mathbb{R}}(\rho_0(x)) \leq \bar{\rho}\left(1 + \tilde{\omega}^{-1}\min\left(0, \inf_{x \in \mathbb{R}}(v_0'(x))\right)\right) \tag{96}$$

$$\rho_0(x) > 0 \text{ for all } x \in \mathbb{R} \tag{97}$$

*Then, the initial-value problem (59), (82), (94) has a unique solution that satisfies the estimates:*

$$\sup_{x \in \mathbb{R}}(\rho(t,x)) \leq \bar{\rho} \text{ , for all } t \geq 0 \tag{98}$$

$$\sup_{x \in \mathbb{R}}\left(|v(t,x) - v^*|\right) \leq \exp(-\tilde{\omega}t)\sup_{x \in \mathbb{R}}\left(|v_0(x) - v^*|\right) \text{ , for all } t \geq 0 \tag{99}$$

$$\rho(t,x) > 0 \text{ , for all } t \geq 0, x \in \Re. \tag{100}$$

*Moreover, there exists a function $f_{TW}: \mathbb{R} \to (0, \bar{\rho}]$ of class $C^1(\mathbb{R}) \cap L^{\infty}(\mathbb{R})$ for which the following estimate holds:*

$$\sup_{x \in \mathbb{R}, t \geq 0}\left(|\rho(t,x) - f_{TW}(x - v^*t)|\exp(\tilde{\omega}t)\right) < +\infty. \tag{101}$$

Notice that the density converges exponentially (in the sup norm) to a traveling wave. This is due to the fact that the whole solution profile ultimately moves with constant speed $v^*$. It should be noted that (96) does not imply that the density is small in general, since $\bar{\rho} \in (0, \rho_{max})$ and can be near $\rho_{max}$. The above convergence result is the macroscopic analogue of the microscopic estimate given by (36) and is not a local result, but a regional result. The estimates given by (98), (99), (100), and (101) are valid for all initial conditions in the region described by (95), (96) and (97). This region is large, since there is no restriction on the speed $v(x)$ and no restriction on how small $\rho_0(x)$ can be.

A significant feature of the macro-PRCC is the fact that it is a system of quasilinear PDEs. Moreover, it is possible to derive some interesting differential equalities for classical solutions of this system of PDEs. When the interval $I(t) \subseteq \mathbb{R}$ is given by the equation $I(t) = (\xi_1(t), \xi_2(t))$ where

$$\dot{\xi}_1(t) = \lim_{x \to (\xi_1(t))^+} (v(t,x)) \text{ and } \dot{\xi}_2(t) = \lim_{x \to (\xi_2(t))^-} (v(t,x)) \text{ for all } t \geq 0 \tag{102}$$

then we can define the functionals

$$I_1(t) = \int_{\xi_1(t)}^{\xi_2(t)} \rho(t,x) dx \tag{103}$$

$$I_2(t) = \int_{\xi_1(t)}^{\xi_2(t)} \rho(t,x) \vartheta(v(t,x)) dx \tag{104}$$

$$I_3(t) = \int_{\xi_1(t)}^{\xi_2(t)} \rho(t,x) \Theta(v(t,x)) dx + \int_{\xi_1(t)}^{\xi_2(t)} \rho(t,x) \Phi\left(\frac{m}{\rho(t,x)}\right) dx \tag{105}$$

where $\vartheta : (0, v_{max}) \to \mathbb{R}$ and $\Theta : (0, v_{max}) \to \mathbb{R}_+$ are given by the formulas

$$\vartheta(v) := \int_{v^*}^{v} \tilde{q}(l) dl \tag{106}$$

$$\Theta(v) = \int_{v^*}^{v} (l - v^*) \tilde{q}(l) dl. \tag{107}$$

The functionals of (103), (104) and (105) are related to the total mass, the total speed and the total mechanical energy, respectively, of the vehicle system. If a classical solution of (58) and (71) satisfies the following conditions for all $t \geq 0$

$$\lim_{x \to (\xi_1(t))^+} (\rho(t,x)) \leq \bar{\rho}, \quad \lim_{x \to (\xi_2(t))^-} (\rho(t,x)) \leq \bar{\rho} \tag{108}$$

then we get (using (58), (71), (102)–(108) and (73), (75)) the following differential equalities for all $t \geq 0$:

$$\dot{I}_1(t) = 0 \tag{109}$$

$$\dot{I}_2(t) = -\int_{\xi_1(t)}^{\xi_2(t)} \rho(t,x) f(v(t,x) - v^*) dx \tag{110}$$

$$\dot{I}_3(t) = -\int_{\xi_1(t)}^{\xi_2(t)} g'(v(t,x))\mu(\rho(t,x))v_x^2(t,x)dx$$
$$-\int_{\xi_1(t)}^{\xi_2(t)} \rho(t,x)\left(v(t,x)-v^*\right)f\left(v(t,x)-v^*\right)dx \tag{111}$$

Notice that since $\mu(\rho) \geq 0$ for $\rho \in (0, \rho_{\max})$, $vf(v) \geq 0$ for $v \in (0, v_{\max})$ (recall (28)) and since $g'(x) > 0$ for all $x \in \mathbb{R}$ (recall (18)), (111) implies that $\dot{I}_3(t) \leq 0$ for all $t \geq 0$. Thus, the total mechanical energy cannot increase, and the total mass remains constant (recall (109)).

If the parameters of the CC that is applied to the vehicles are selected in such a way that
$$P'(\rho) \geq 0, \text{ for all } \rho \in (0, \rho_{\max}) \tag{112}$$
$$g(v) = v, \text{ for all } v \in (0, v_{\max}) \tag{113}$$
$$\mu(\rho) = \frac{1}{\tilde{k}}\rho P'(\rho) \tag{114}$$
$$f(v-v^*) = \tilde{k}\vartheta(v), \text{ for all } v \in (0, v_{\max}) \tag{115}$$

where $\tilde{k} > 0$ is a constant, then (104) and (110) imply that, if (108) holds, then
$$\dot{I}_2(t) = -\tilde{k}\, I_2(t) \Rightarrow I_2(t) = \exp\left(-\tilde{k}\, t\right)I_2(0), \text{ for } t \geq 0. \tag{116}$$

Therefore, in this case (where (108), (112) – (115) hold) the total speed converges exponentially to zero. The particular case where (58), (71), (112), (113), (114), and (115) hold was studied in [95]. If we define
$$\bar{\varphi}(t,x) = \vartheta(v(t,x)) + \rho^{-2}(t,x)\mu(\rho(t,x))\rho_x(t,x), \text{ for } t \geq 0, \ x \in I(t) = (\xi_1(t), \xi_2(t)) \tag{117}$$

then, by virtue of (58), (71), (113), (114) and (115), we obtain the differential equation
$$\bar{\varphi}_t + v\bar{\varphi}_x = -\tilde{k}\bar{\varphi}. \tag{118}$$

The similarity of (118) with (65) is striking. Moreover, in this case we can define the additional functional
$$I_4(t) = \frac{1}{2}\int_{\xi_1(t)}^{\xi_2(t)} \rho(t,x)\bar{\varphi}^2(t,x)dx \tag{119}$$

And, using (58), (102) and (118), we show that
$$\dot{I}_4(t) = -2\tilde{k}I_4(t) \Rightarrow I_4(t) = \exp\left(-2\tilde{k}\, t\right)I_4(0), \text{ for } t \geq 0. \tag{120}$$

It becomes clear that the subset of the state space where $\bar{\varphi}(t,x) \equiv 0$ is invariant and exponentially attracting. On this invariant set, the solutions satisfy the equation
$$\rho_t + \left(\rho\vartheta^{-1}\left(-\frac{P'(\rho)}{\tilde{k}\rho}\rho_x\right)\right)_x = 0 \tag{121}$$

where $\vartheta^{-1}: \mathbb{R} \to (0, v_{\max})$ is the inverse function of $\vartheta: (0, v_{\max}) \to \mathbb{R}$ defined by (106) (well-defined because $\vartheta((0, v_{\max})) = \mathbb{R}$). The nonlinear PDE of (121) is not hyperbolic and is not parabolic. When $\rho(t,x) \leq \bar{\rho}$, we get that $\rho_t + v^*\rho_x = 0$ (a transport PDE) with propagation speed being equal to $v^*$; when $\rho(t,x) > \bar{\rho}$, we get a nonlinear convection-diffusion equation with the diffusion coefficient depending on both the density $\rho$ and its spatial derivative $\rho_x$. This feature is rarely studied in the literature, where the diffusion coefficient exclusively depends on the density $\rho$, see for instance [79], [85], [98], with notable exceptions are [62], [63], [72], where the viscosity also depends on the spatial derivative of the state. The solutions of (121) approximate the solutions of the macroscopic model given by (58), (71), (112) – (115). Lyapunov functionals and a finite-difference numerical scheme for the PDE given by (121) were provided in [95].

The analogy with traffic flow of human-driven vehicles is now clear. For human-driven vehicles, the ARZ model has an invariant and exponentially attracting set, on which the traffic flow is described by the LWR model. For CAVs, the macro-PRCCs with the selections given by (112) – (115) has an invariant and exponentially attracting set. On this invariant set, the traffic flow is described by (121). Therefore, (121) is a reduced model for CAVs; loosely speaking, it is the "LWR model for CAVs".

## 4. Simulations

### 4.1 Viscous and Inviscid NCC

Consider $n=15$ vehicles of length $\sigma_i = 5m$ with desired speed $v^* = 30 m/s$ on a road of width $2a = 14.4m$ and speed limit $v_{\max} = 35 m/s$. Let $L = 5.59m$, $\lambda = 25m$, $\varphi = 0.25$, $\gamma = 0.1$ and $\Gamma = 0.5$. The constants $L$ and $p_{i,j}$ were selected by following the optimal selection of safety distance and eccentricity described in [69]. We also select

$$V_{i,j}(d) = \begin{cases} z_1 \dfrac{(\lambda - d)^3}{d - L} & , L < d \leq \lambda \\ 0 & , d > \lambda \end{cases} \tag{122}$$

$$U_i(y) = \begin{cases} \left( \dfrac{1}{a^2 - y^2} - \dfrac{\bar{c}}{a^2} \right)^4 & , -a < y < -\dfrac{a\sqrt{\bar{c}-1}}{\sqrt{\bar{c}}} \text{ and } \dfrac{a\sqrt{\bar{c}-1}}{\sqrt{\bar{c}}} < y < a \\ 0 & -\dfrac{a\sqrt{\bar{c}-1}}{\sqrt{\bar{c}}} \leq y \leq \dfrac{a\sqrt{\bar{c}-1}}{\sqrt{\bar{c}}} \end{cases} \tag{123}$$

with $\bar{c} = 1.5$ and $z_1 = 10^{-4}$. For the viscous case, we let $g(s) = s$ and

$$\kappa_{i,j}(d) = \begin{cases} z_2 (\lambda - d)^2 & L < d \leq \lambda \\ 0 & d > \lambda \end{cases}, \text{ for all } i, j = 1, \ldots, n \tag{124}$$

with $z_2 = 0.03$. Finally, the function $\ell(x)$ satisfying (17) is given by

$$\ell(x) = \dfrac{1}{2\eta} \begin{cases} 0 & \text{if } x \leq -\eta \\ (x + \eta)^2 & \text{if } -\eta < x < 0 \\ \eta^2 + 2\eta x & \text{if } x \geq 0 \end{cases} \tag{125}$$

with $\eta = 0.2$.

Fig. 4 shows the norm $|(v_1(t) - v^*, \ldots, v_{15}(t) - v^*)|_\infty$ for both the viscous and inviscid NCC illustrating the convergence of the speeds of the vehicles to the desired speed $v^*$. Fig. 5 shows the maximum value of the acceleration over time among all vehicles (i.e., $\max_{i=1,\ldots,15}(|F_i(t)|)$) and its convergence for both the inviscid and viscous cases of the NCC. Finally, Fig. 6 shows the minimum inter-vehicle distance, $\min_{i,j=1,\ldots,n}(d_{i,j}(t))$, indicating that there are no collisions among vehicles, since $\min_{i,j=1,\ldots,n}(d_{i,j}(t)) > L$. To obtain the numerical solution of the closed-loop system of $n$ vehicles described by (1) with the NCC (and also for PRCC and GCC in the next sections), we have used an adaptive step-size numerical scheme based on Euler and Heun method (see [35]) with initial conditions selected from the set $\Omega$ in (3) with $v_{\max}$, $\varphi$, and $a$ as given above. For better appreciation, the video in (https://youtu.be/eTtJT_k42q4) illustrates the two-dimensional movement of vehicles on lane-free roads using the viscous NCC. The video shows that vehicles do not collide with each other, they remain within the road boundaries, and eventually reach the same desired speed.

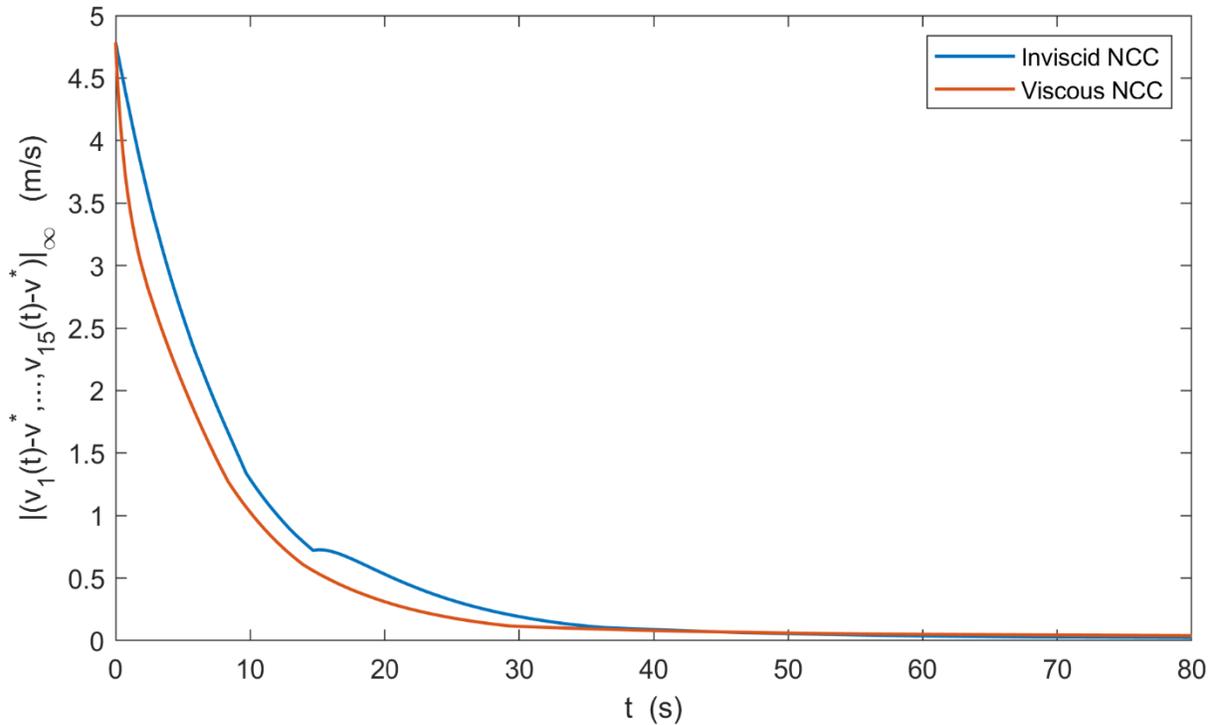

**Fig. 4.** Convergence of the speeds to the desired speed for the viscous and inviscid NCC

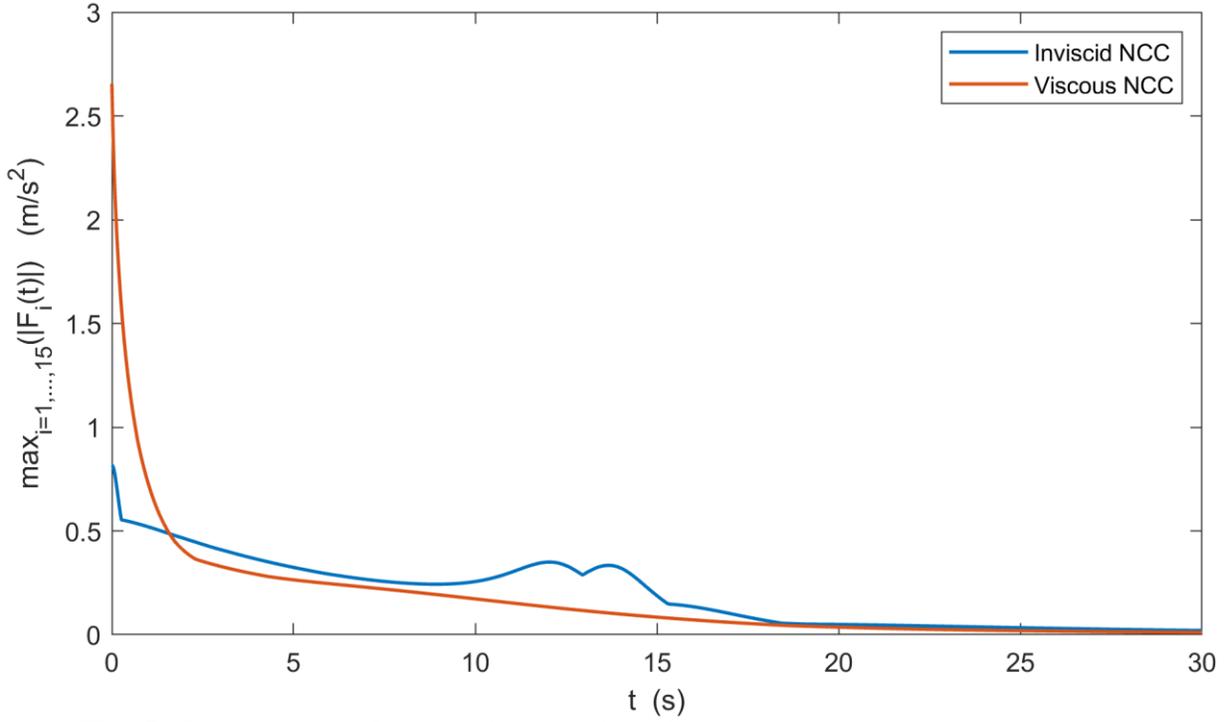
**Fig. 5:** Convergence of the maximum value of the acceleration among all vehicles.

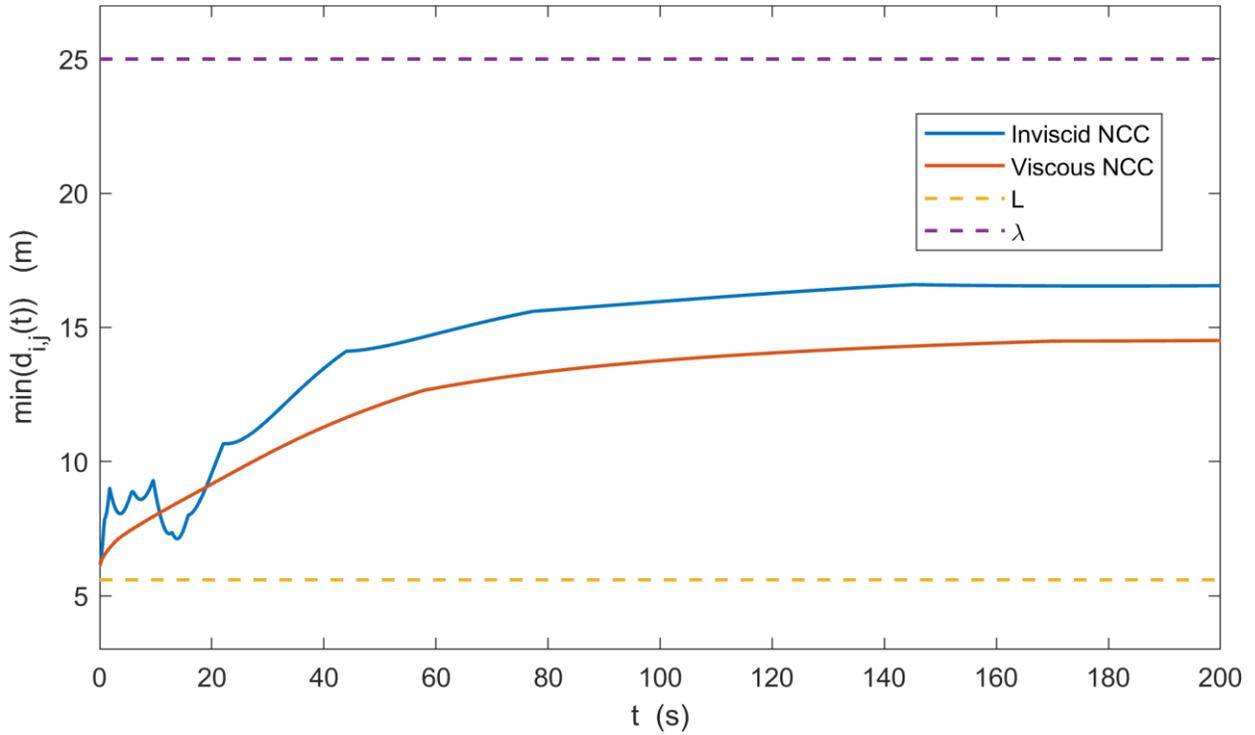
**Fig. 6:** Minimum inter-vehicle distance.

## 4.2 Viscous and Inviscid PRCC

For the PRCC, we consider the case where $f(s)=0.5s$, $\bar{f}(s)=2s$, and $A=1$, $b=1$. We also consider the potential functions in (122) and (123) with $z_1=0.01$ and $\bar{c}=1.5$. For the viscous PRCC we consider $\kappa_{i,j}$ as defined in (124) with $z_2=0.1$. Fig. 7 shows the norm $|(v_1(t)-v^*,...,v_{15}(t)-v^*)|_\infty$ for both the viscous and inviscid PRCC illustrating the convergence of the speeds of the vehicles to the desired speed $v^*$. Fig. 8 shows the maximum value of the acceleration $F_i$ among all vehicles over time as well as its convergence to zero. Finally, Fig. 9 shows the minimum inter-vehicle distance verifying that there are no collisions among vehicles for both the viscous and inviscid PRCC. The video (https://youtu.be/3qrWxRt7xKQ) illustrates the qualitative properties of the viscous PRCC, where vehicles do not collide with each other, remain within the road boundaries, and eventually reach the same desired speed.

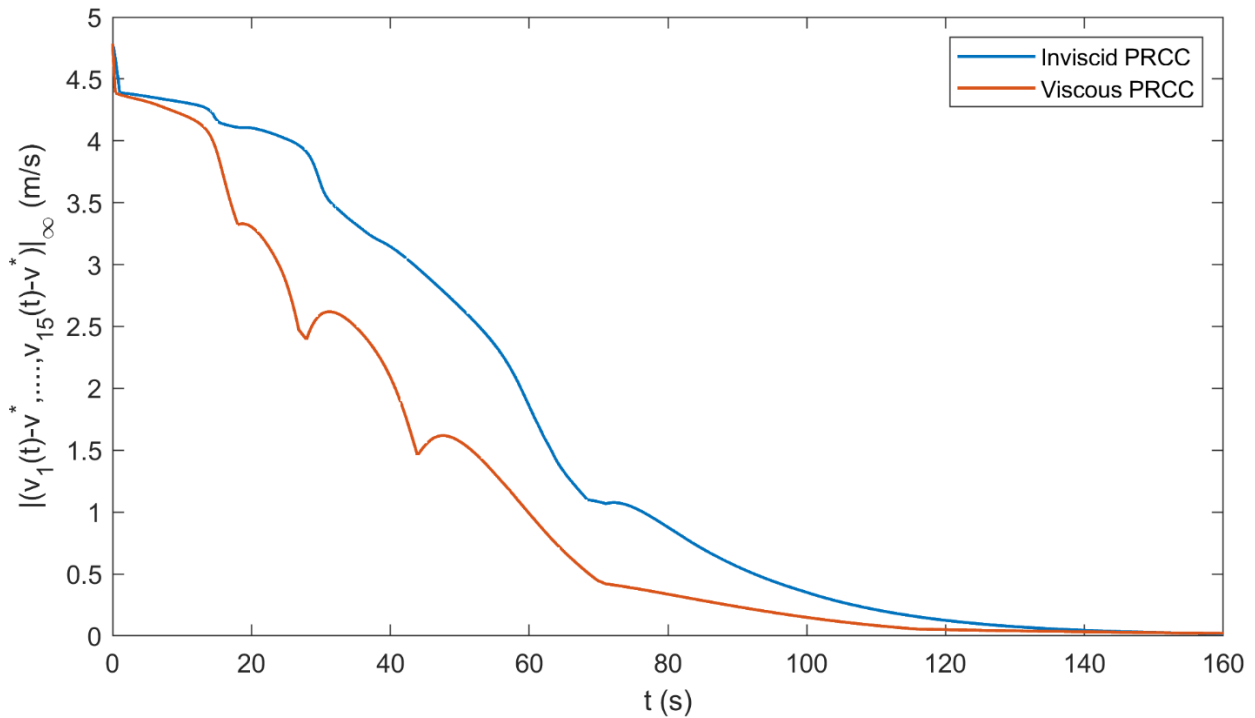

**Fig. 7:** Norm of speeds $|(v_1(t)-v^*,...,v_{15}(t)-v^*)|_\infty$ for the viscous and inviscid PRCC.

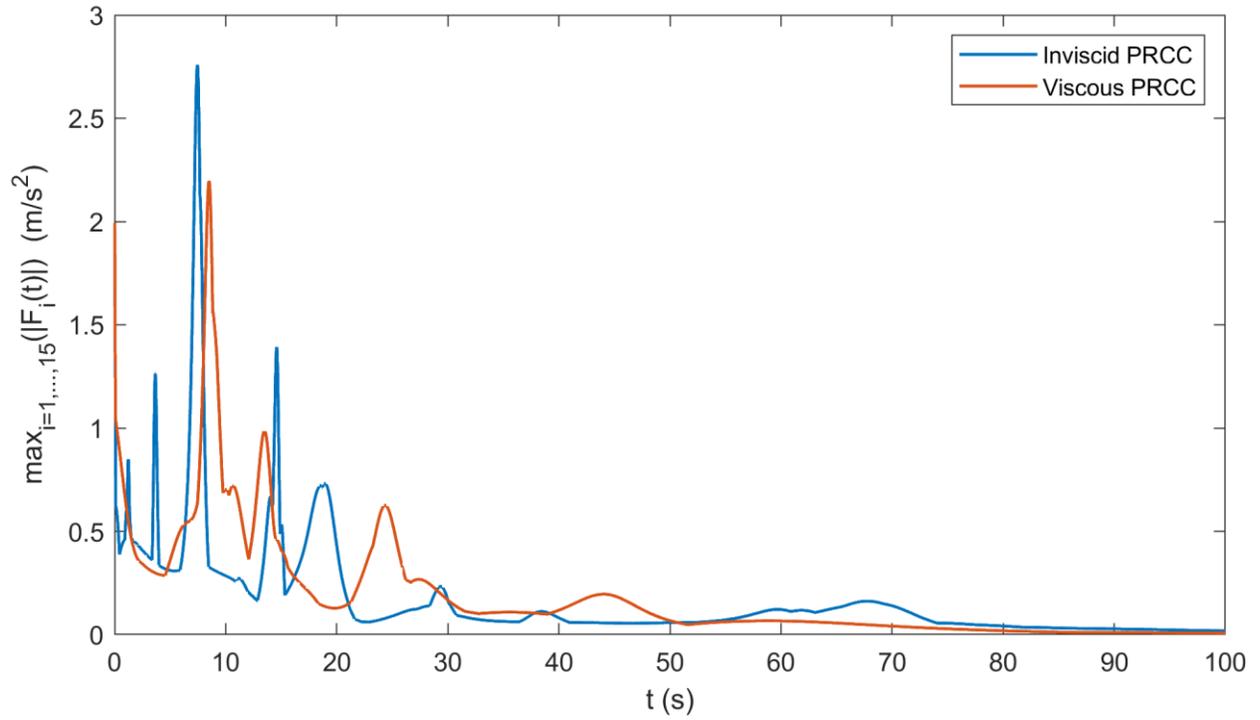

**Fig. 8:** Maximum value of acceleration among all vehicles over time.

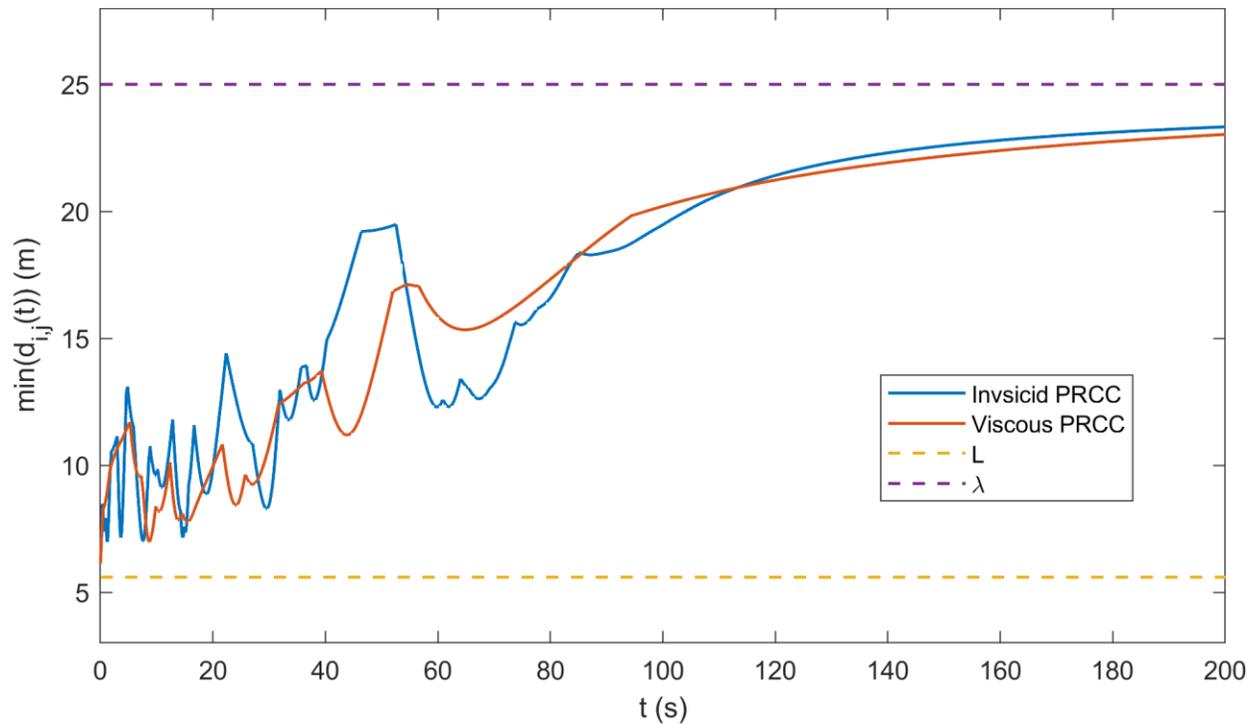

**Fig. 9:** Minimum inter-vehicle distance.

## 4.3 Application of the GCC on a Road with an Off-ramp

To illustrate the application of the GCC, we consider $n = 150$ vehicles of the same length $\sigma_i = 5m$ on a lane-free road with an off-ramp. The corridor boundaries $\alpha_i(x)$ and $\beta_i(x)$ were appropriately selected for certain vehicles to exit the highway. Moreover, the speed set-point $v_i^*$ for each vehicle is selected randomly from the set $\{28, 29, 30, 31\}$ in m/s. Fig. 10 shows the road layout with the off-ramp as well as simulation snapshots at different time instances, whereby the vehicles with black color are bound to exit the highway by following the corridor highlighted by red-dashed lines. The video (https://youtu.be/VAbvXBtCuaA) shows the complete simulation of the vehicles exiting the highway. The minimum inter-vehicle distance is shown in Fig. 11, illustrating that there are no collisions between vehicles since $\min_{i,j=1,\ldots,n} (d_{i,j}(t)) > L$. For the application of the GCC, the following parameters were used $\gamma = 0.1$, $\Gamma = 1$, $b = 1$, $\varepsilon = 0.001$, $p_{i,j} = 4.25$, $i, j = 1, \ldots, n$, $\varphi = 0.45$, $v_{\max} = 35 m/s$, and $c = 2.1$, $z_1 = 10^{-4}$, $\lambda = 100 m$, $L = 6m$ for the potentials defined in 122 and 123. An extra video (https://youtu.be/ZdsSxe1Zfa8) shows the application of the GCC in a more complex scenario with a road that includes an off-ramp, an on-ramp, and has variable width. The video shows that vehicles with black color follow the corridor with red-dashed lines towards the off-ramp and re-enter the highway via the on-ramp, while the rest of the vehicles adjust their position to remain in the highway when its width is reduced.

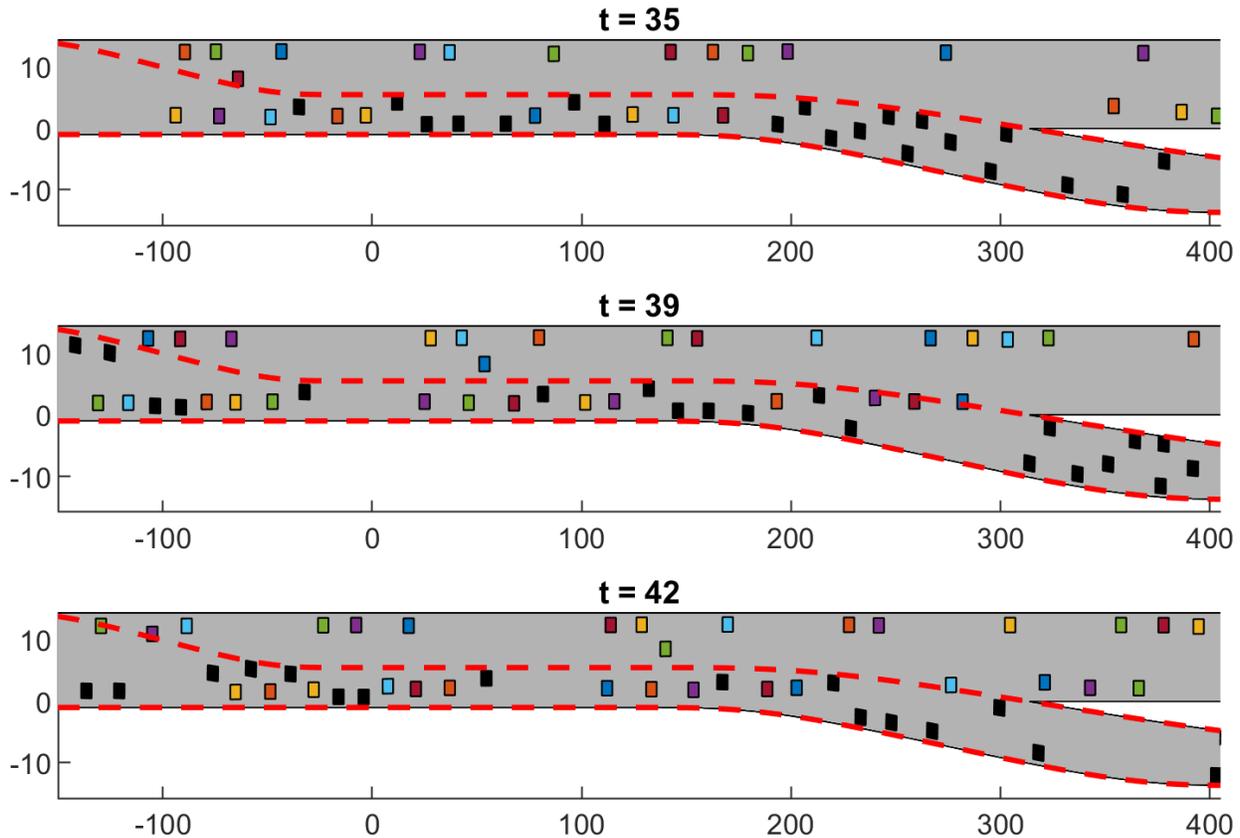

**Fig. 10:** Simulation outcome at various time instances.

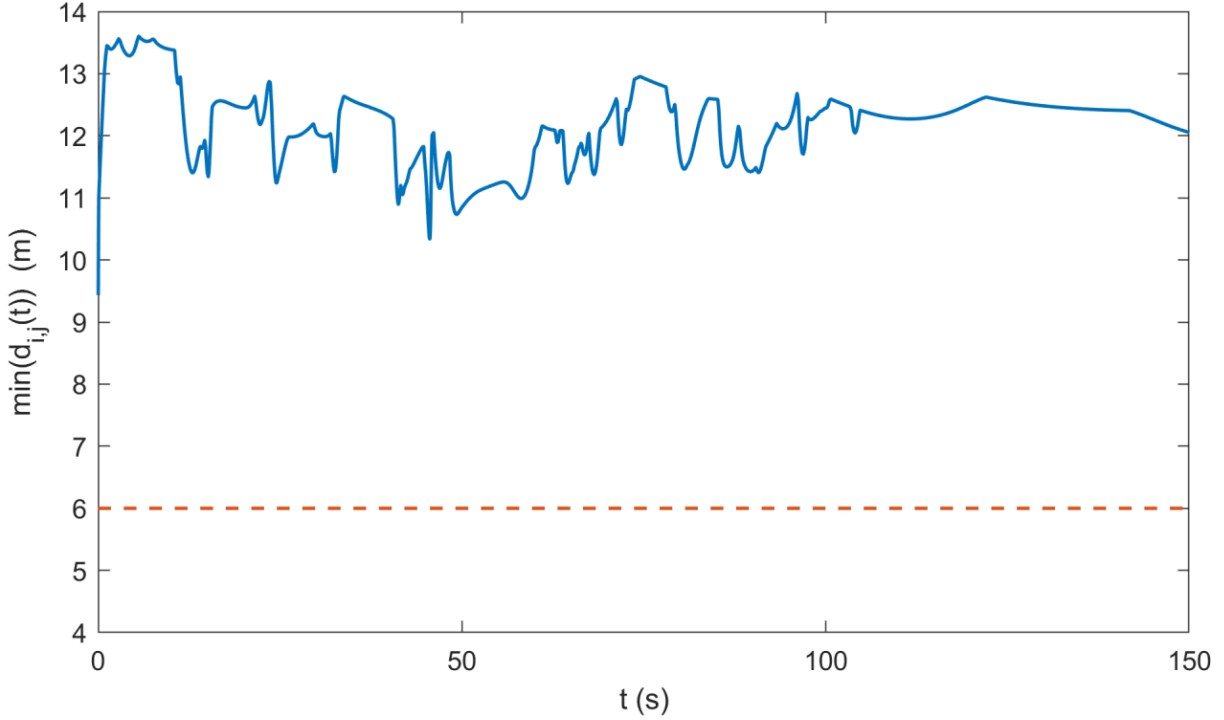

**Fig. 11:** Inter-vehicle distance for the application of the GCC.

### 4.4. Topological Obstructions of the GCC

We consider two vehicles with $x_i(0) = 0$, $v_i(0) = 30m/s$, $\theta_i(0) = 0$, $i = 1, 2$, and $y_1(0) = 5$, $y_2(0) = -5$ in a road section with $\beta(x) = -\alpha(x)$ and $\beta(x) - \alpha(x) = 4$, for $x > 80$. Moreover, we select $p_{i,j} = 1$, $i, j = 1, 2$, $L = 6m$, and $v^* = 30m/s$. Notice that both the road boundaries and the positions of the vehicles are completely symmetric with respect to $y = 0$, and for $x > 80$ the width of the road is constant and small enough so that only one vehicle can be placed laterally. This scenario illustrates that in very specific cases, the GCC does not guarantee a lower bound on the speeds of vehicles. However, by slightly changing the initial positions of the vehicles, both vehicles can pass through the bottleneck, see Fig. 12.

### 4.5 Nonlinear Heat Equation

We consider a traffic scenario and compare the density and flow of (121) with one of the most well-known traffic flow models for human drivers, the LWR model given by (58) and (59) with

$$F(\rho) = v_f \exp\left(-\frac{1}{\hat{a}}\left(\frac{\rho}{\rho_c}\right)^{\hat{a}}\right) \qquad (126)$$

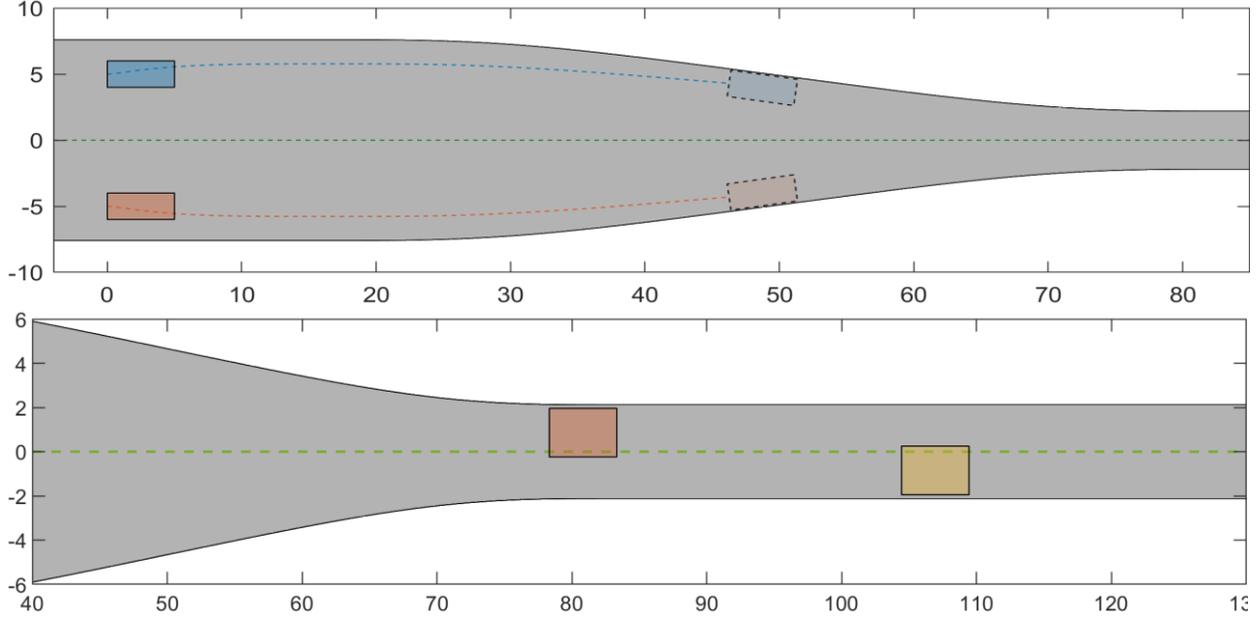

**Fig. 12:** Vehicles coming to a halt for specific initial conditions (top) and vehicles passing through the bottleneck for perturbed initial conditions (bottom).

where $v_f > 0$ is the free speed, $\rho_c$ is the critical density and $\hat{a} > 0$ is a parameter. We consider a single-lane motorway with $v_f = 102 \ km/h$, $\rho_{max} = 180 \ veh/km$, $\rho_c = 33.3 \ veh/km$, and $\hat{a} = 2.34$ that were estimated in [60] based on real data from a part of the Amsterdam A10 motorway. In this scenario, we assume that the initial density on the road is characterized by a congestion belt in the interval $[1.5, 2.75]$ (km) as shown in Fig. 13 with $\rho_0(x) = 0$ for $x < 0$ and $x > 4$. For (121) we select $\bar{\rho} = 31 \ veh/km$, $k = 1/40$, $v^* = 102 \ km/h$ and

$$P'(\rho) = \begin{cases} 0 & , 0 \le \rho \le \bar{\rho} \\ \dfrac{(\rho - \bar{\rho})^2}{\rho_{max} - \rho} & , \bar{\rho} < \rho < \rho_{max} \end{cases}$$

For the LWR model, we consider the interval $x \in [0,120]$ and use Godunov's method to obtain its solution. The density profiles for both models over a simulation time of $1h$ are shown in Figs. S13 and S14, respectively, where the density of both models dissipates along the road. However, for the LWR model the dissipation is stronger and the density is spread over a large road interval for increasing time. More specifically, for $t = 1 \ (h)$, the density is non-zero over the interval $[90,106]$ (km), which implies that the vehicles retain large inter-vehicle distances while the maximum density is $\max_{x \in [90,106]} (\rho[1]) = 12.3 \ (veh/km)$. For the model (121), the density dissipates at a lower rate but the vehicles remain in a $4km$ stretch of the road ($x \in [102,106]$ (km), since the desired speed is $v^* = 102 \ (km/h)$), as was the case in the initial density, and with maximum density $\max_{x \in [102,106]} (\rho[1]) = 37.19 \ (veh/km)$ (which eventually approaches the critical density $\rho_c$ since $\bar{\rho} < \rho_c$).

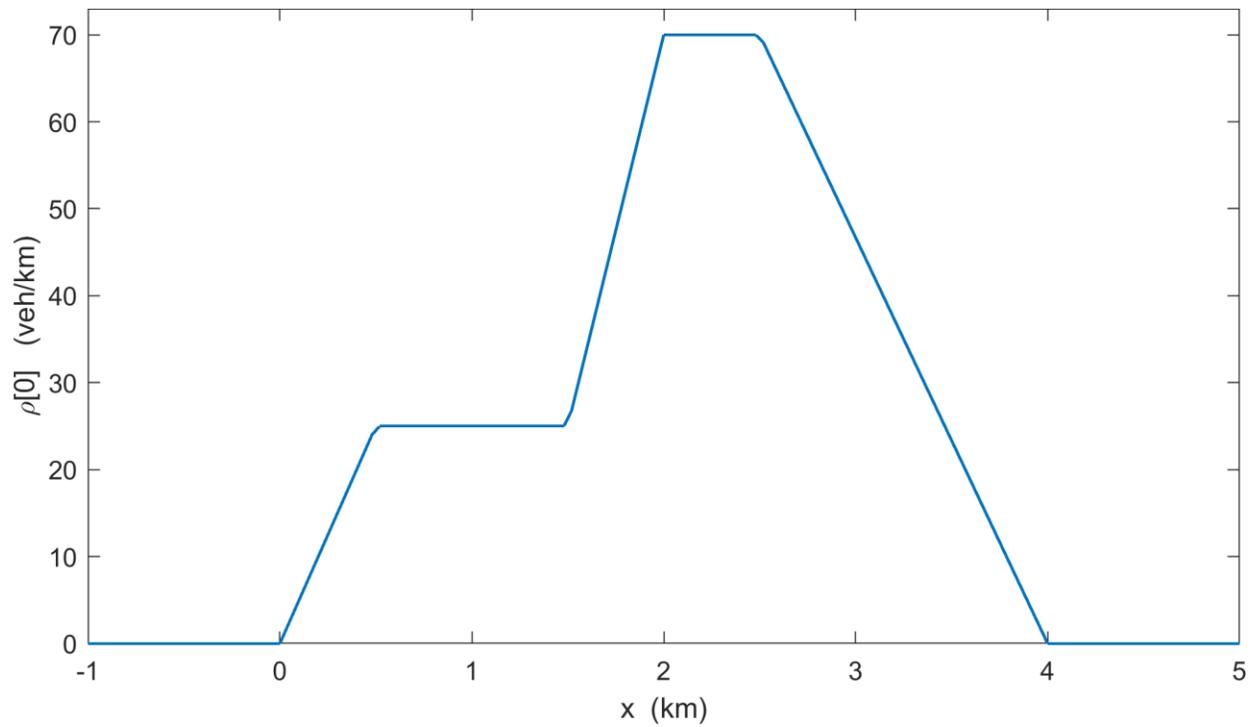
**Fig. 13:** Initial density $\rho[0]$

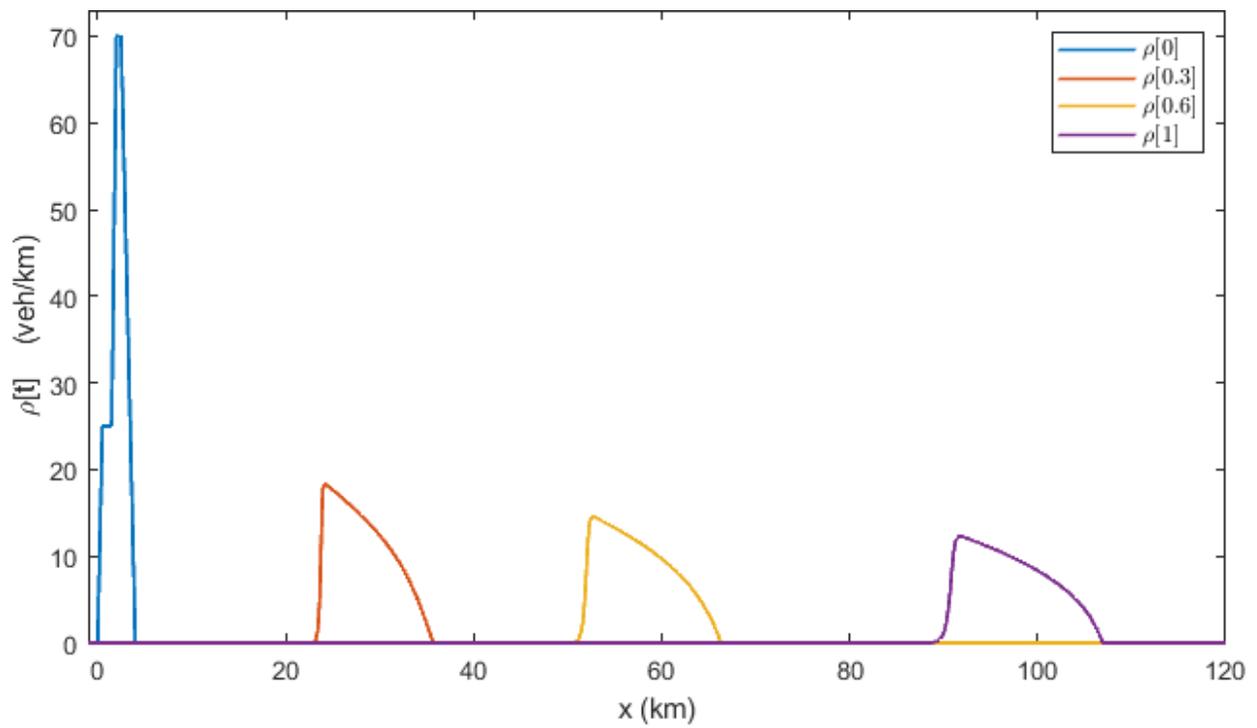
**Fig. 14:** Density profiles for the LWR.

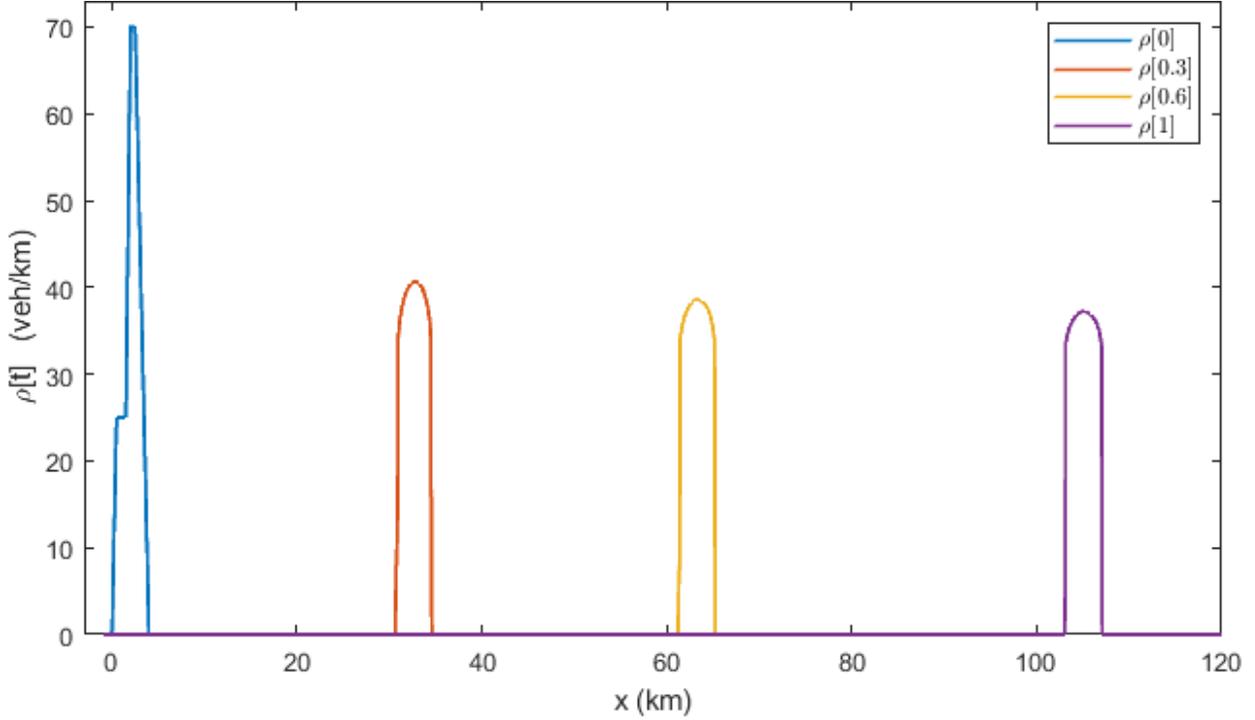

**Fig. 15:** Density profiles for model (121).

Since with the model (121) the density converges towards the critical density while vehicles are spread in a small space interval, it is expected that the flow $\rho v$ of the model (121) is much higher than the flow of the LWR model. Equivalently, vehicles equipped with the PRCC retain smaller inter-vehicle distances (higher density) than in the case of human-driven vehicles. Indeed, to compare the mean flow for each model, we define

$$\text{Mean Flow } = \frac{1}{T}\int_0^T \int_{m_S(t)}^{m_F(t)} \frac{\rho(t,x)v(t,x)}{m_F(t)-m_S(t)} dx dt$$

where $T > 0$ denotes the time horizon and $[m_S(t), m_F(t)]$ denotes the interval where $\rho(t,x) \neq 0$ for $x \in [m_S(t), m_F(t)]$ at each time instant $t \geq 0$. The mean flow for both models and $T = 1$ ($h$) is shown in Table 3. The mean flow for the model (121) is much higher than the mean flow of the LWR, even when the desired speed is much lower ($v^* = 51 km/h$). The density profiles of the CAV model with $v^* = 102 km/h$ are shown in Fig. 15. Note that when vehicles are moving with free-flow speed, the overall travel time is lower with the CAV model, since for $t = 1$ (h), vehicles are in the interval $x \in (102, 106)$ (km) whereas with the LWR, vehicles are in the interval $x \in (90, 106)$ (km).

**Table 3. Comparison of mean flows.**

|  | Model (122) with $v^* = 51 km/h$ | Model (122) with $v^* = 102 km/h$ | LWR with $v_f = 102\ km/h$ |
|---|---|---|---|
| Mean Flow | 1343 veh/h | 2833 veh/h | 965 veh/h |

## 4.6 Traveling Waves of the Inviscid macro-NCC

The following example illustrates certain properties of the inviscid macro-NCC, which show that the speed converges exponentially to the speed set-point $v^*$; while the density converges to a travelling wave. We consider a road with initial density and initial speed given by

$$\rho_0(x) = 25 + \begin{cases} (10x)^2(x-1)^2, & x \in (0,1) \\ 0 & else \end{cases}$$

$$v_0(x) = 60 + \begin{cases} (4x)^3(x-1)^3, & x \in (0,1) \\ 0 & else \end{cases}.$$

These initial conditions indicate that there is a congestion belt on the interval $x \in (0,1)$ where vehicles are moving at lower speed and accelerate to a speed of $v^* = 60 km/h$ as the density decreases to a constant value. Values of $v^* = 60 km/h$ and $\omega = 5$ were used. Fig. 16 displays the speed profiles $v[t]$ (left) and the density profiles $\rho[t]$ (right) at different time instants.

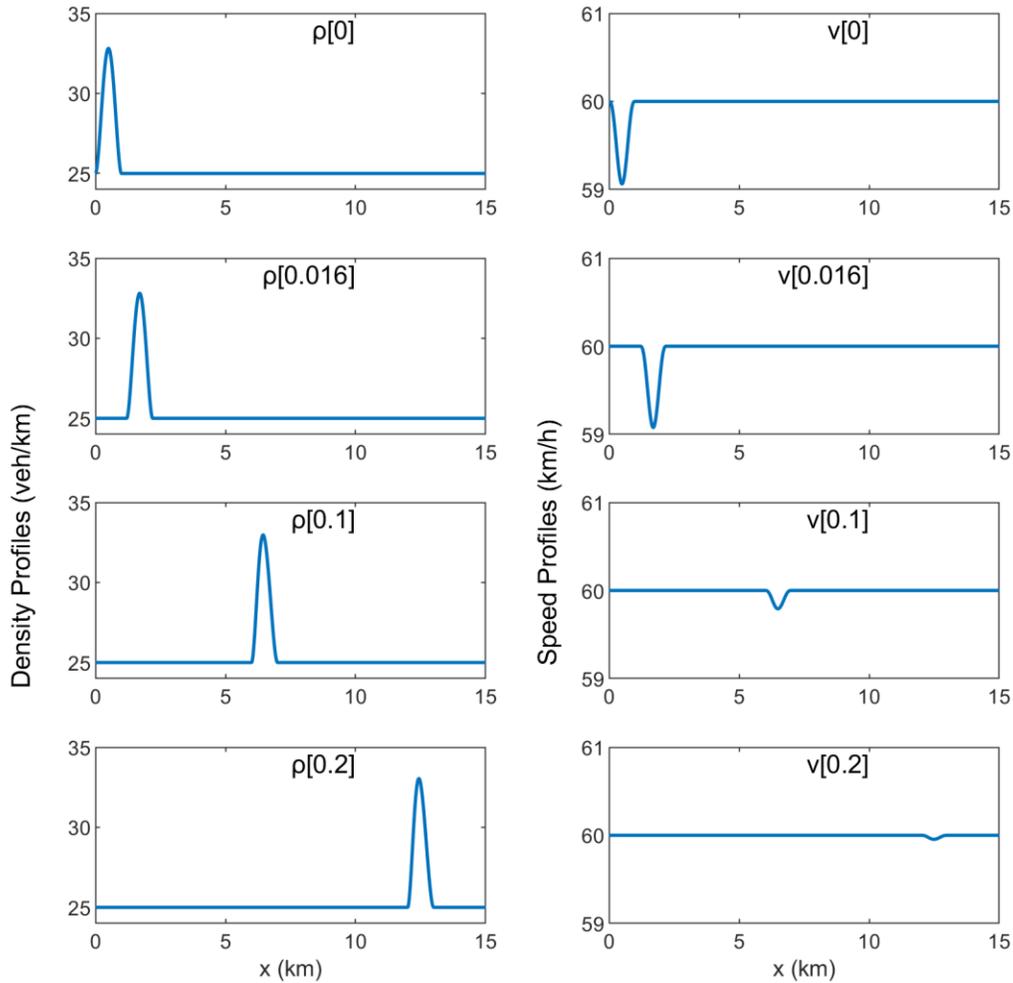

**Fig. 16:** Density profiles (left) and speed profiles (right) of the inviscid macro-NCC.

## 4.7 Macroscopic Model for the PRCCs

In this simulation scenario we compare the macro-PRCC, given by (58), (71), (112) – (115), versus the ARZ model for human-driven vehicles given by (58) and (61). For the macro-PRCC, we select the following function for the pressure

$$P(\rho) = \bar{\rho}(v^*)^2 \begin{cases} 0 & , 0 \leq \rho \leq \bar{\rho} \\ \tilde{k}\left(1 - \dfrac{\rho}{\bar{\rho}} - \dfrac{\bar{\rho}(\rho_{max} - \bar{\rho})^2}{\rho_{max}} \ln\left(\dfrac{\rho_{max} - \rho}{\rho_{max} - \bar{\rho}}\right) + \dfrac{\bar{\rho}}{\rho_{max}} \ln\left(\dfrac{\rho}{\bar{\rho}}\right)\right) & , \bar{\rho} < \rho < \rho_{max} \end{cases}.$$

The numerical solution of the macro-PRCC is obtained using the particle method given by (66) – (70), the corresponding functions from Table 1, (106), (113) – (115), and the following parameters: $\tilde{k} = 30$, $v^* = 33 km/h$, $v_{max} = 35 km/h$, $\bar{\rho} = 63.158 veh/km$, $\rho_{max} = 120 veh/km$. The ARZ model is used with

$$F(\rho) = v_f \exp\left(-\dfrac{1}{\hat{a}}\left(\dfrac{\rho}{\rho_c}\right)^{\hat{a}}\right)$$

and $\bar{k} = 0.05$, $\hat{a} = 2.34$, where $v_f = 33 km/h$ is the free speed and $\rho_c = 63.158 veh/km$ is the critical density. To obtain the numerical solution of the ARZ we used a modified Nessyahu-- Tadmor scheme, [65].

The density and speed profiles of the ARZ and the macro-PRCC are shown in Figs. 17 - 20. For the ARZ, the density is spread on a larger space interval, while the speed decreases at a higher rate at larger densities. In contrast, with the macro-PRCC, the density dissipates but remains in a much smaller space interval while there is no abrupt change in the speed.

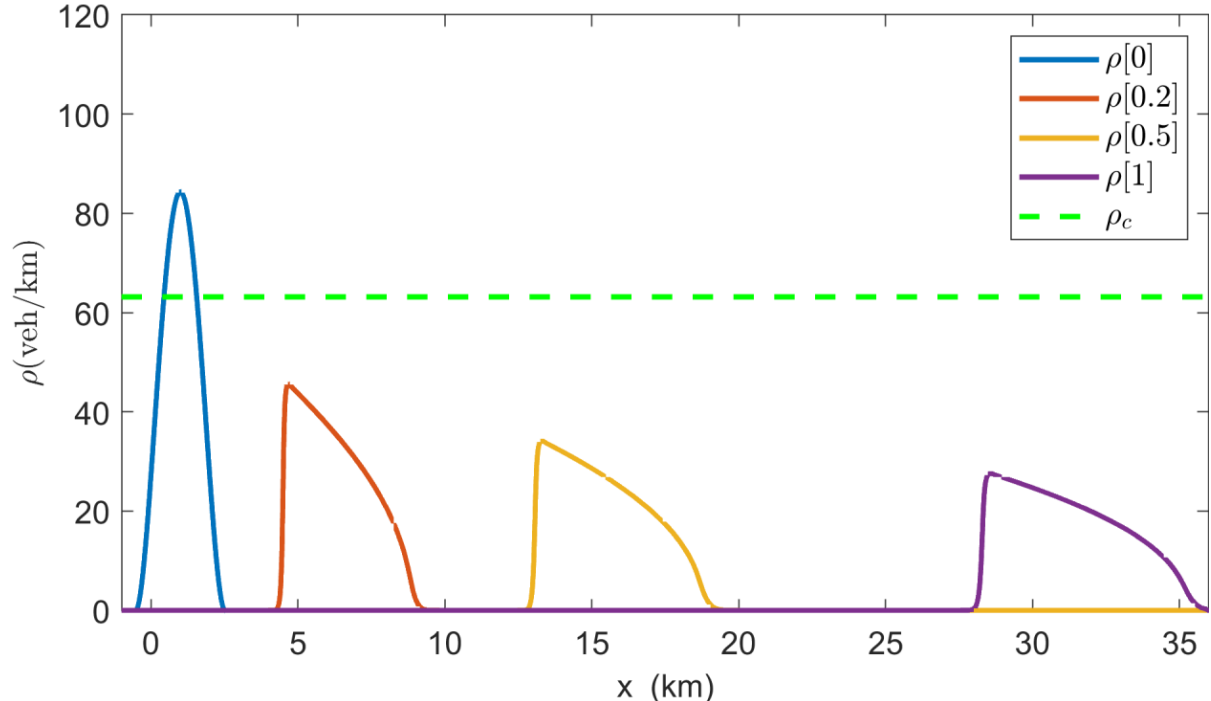

**Fig. 17:** Density profiles $\rho[t]$ for the ARZ.

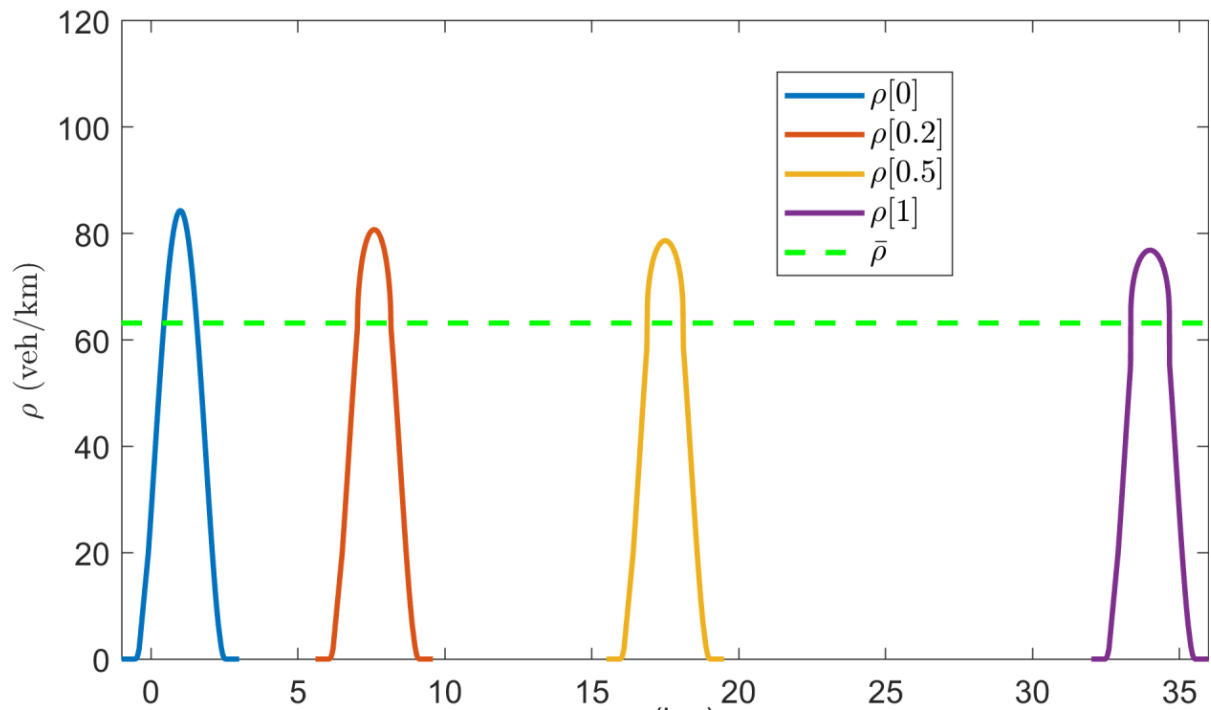

**Fig. 18:** Density profiles $\rho[t]$ for the macro-PRCC.

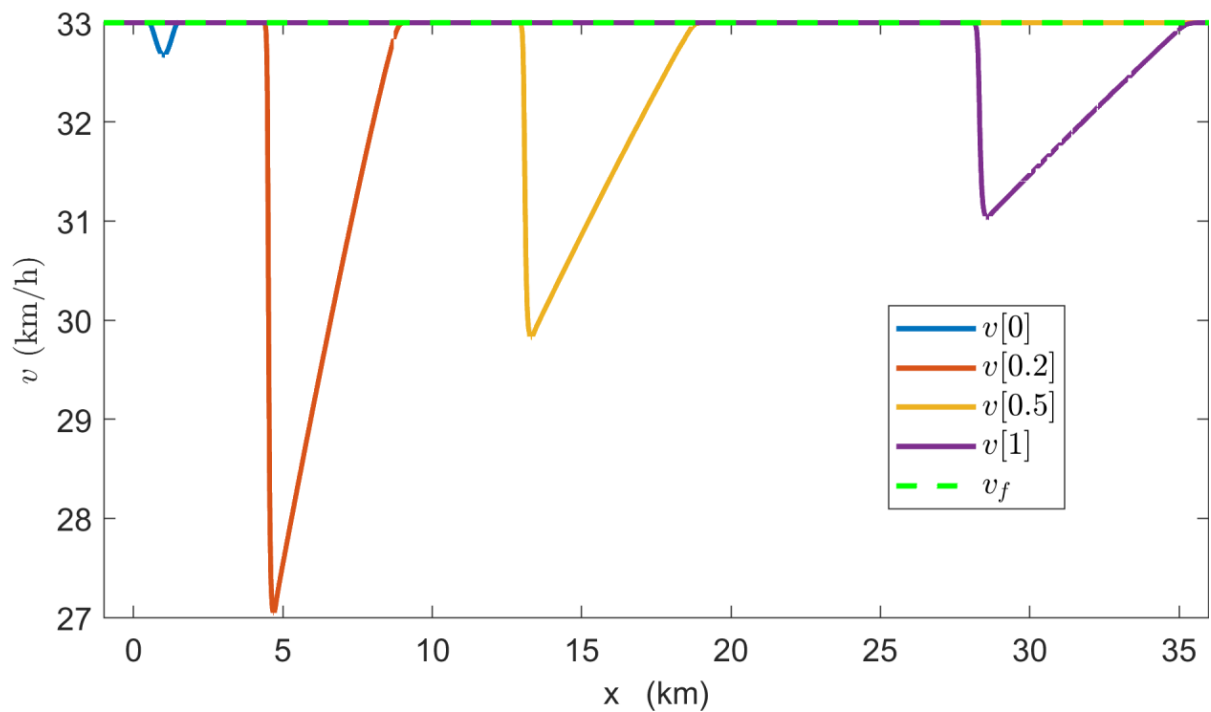

**Fig. 19:** Speed profiles $v[t]$ for the ARZ

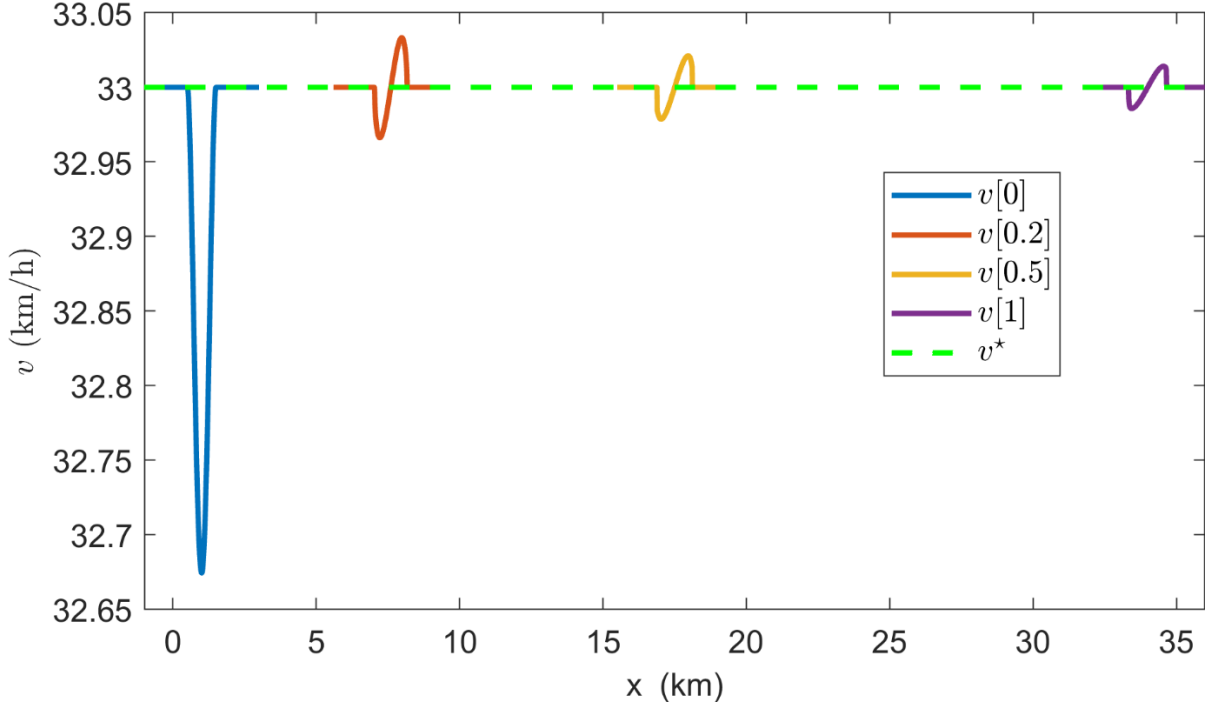

**Fig. 20:** Speed profiles $v[t]$ for the macro-PRCC

## 5. Conclusions and Future Directions

This paper presented in detail families of controllers for the two-dimensional movement of CAVs via algorithmic creation of artificial forces that induce the vehicle motion in lane-free traffic with vehicle nudging, rigorously guaranteeing a variety of desired features, including safety in terms of collision avoidance and road boundary respect. The design of each of the families of CCs is based on a control Lyapunov methodology with the Lyapunov function expressed on measures of the total energy of the system. The proposed families of cruise controllers are decentralized (per vehicle) and require either the measurement only of the relative distances from adjacent vehicles (inviscid cruise controllers) or the measurement of relative speeds and relative distances from adjacent vehicles (viscous cruise controllers). Moreover, we discussed the formal derivation of the corresponding macroscopic models from the derived families of CCs by using particle methods and provided direct relations between selectable CC features and the resulting macroscopic traffic flow characteristics. This allows the active design of efficient traffic flow with desired properties, i.e., the construction of artificial traffic fluids.

The presented families of CCs (NCC, PRCC, GCC) contain a wide range of selectable functions and parameters that may lead to very different behaviors at the microscopic and macroscopic levels. This calls for more detailed exploration of the offered options based on appropriate criteria and secondary goals. Furthermore, radically novel CCs, offering similar qualitative guarantees as those presented here, might provide additional benefits, [86], [9].

The general setting and derivation of the CCs presented earlier, is not limited to automobile lane-free traffic, but can be further extended for application to a variety of moving objects, ranging from 3-D lane-free movement of aerial vehicles (drones) to 2-D lane-free vessel movement.

While vehicle automation has made tremendous advances in the last decade, with a variety of driving-assistance systems already available in today's market, the path to fully autonomous vehicles requires first the co-existence of human-driven and "self-driven" vehicles. Although the CCs presented in this review, were designed for fully automated vehicles with lane-free roads in mind, their decentralized nature can also handle cases of mixed traffic. Thus, the applicability of the CCs needs to be further strengthened by rigorous mathematical guarantees, not only for the not-too-distant future of automated vehicles, but also for the intermediary "mixed" traffic conditions, [44], [46], [108].

There are still challenges in the deployment of CAVs in urban road networks. While autonomous vehicles can use decentralized CCs for their movement, they are not entirely isolated entities driving in a fully protected environment. Uncertainties in a dynamic traffic environment, such as pedestrians that have free and unpredictable movements, require important advancements in both technology and pedestrian motion prediction.

Finally, the rigorous (rather than formal) derivation of the macroscopic model for lane-free CAV traffic with CCs remains an open problem. This is expected since the rigorous derivation of the NS equations from the microscopic motion of fluid particles is also an open problem (Hilbert's 6th problem, [42]).

**Acknowledgments**

The research leading to these results has received funding from the European Research Council under the European Union's Horizon 2020 Research and Innovation programme/ ERC Grant Agreement n. [833915], project TrafficFluid.

Table 2: Macroscopic Models for CAVs.

| | macro-PRCC | macro-NCC | macro-GCC | 1-D NS for a polytropic gas |
|---|---|---|---|---|
| Continuity | $\rho_t + (\rho v)_x = 0$ | | | |
| Speed | $Q(v,r)(v_t + v v_x) = G - J f(v - V(r))$, $r = \rho^{-1} P'(\rho) \rho_x$, $G = \rho^{-1} S(v,r)\left(\mu(\rho) g'(v) v_x\right)_x - r$ | | | |
| States | $\rho \in (0, \rho_{\max})$, $v \in (0, v_{\max})$ | | | $\rho > 0$, $v \in \mathbb{R}$ |
| Pressure | $P:(0,\rho_{\max}) \to \mathbb{R}$, $\int_{\bar{\rho}}^{\rho_{\max}} \rho^{-2} P(\rho) d\rho = +\infty$, $P(\rho) = 0$ for $\rho \in (0, \bar{\rho}]$ | | | $P(\rho) = B \rho^\eta$ |
| Constants | $\rho_{\max}, v_{\max}, \gamma, \varepsilon, M > 0$, $v^* \in (0, v_{\max})$, $\bar{\rho} \in (0, \rho_{\max})$ | | | $B > 0$, $\eta \in (1,2)$ |
| Viscosity | $\mu:(0, \rho_{\max}) \to \mathbb{R}_+$, $\mu(\rho) = 0$ for $\rho \in (0, \bar{\rho}]$ | | $\mu(\rho) = \rho P'(\rho)$ | $\mu:(0,+\infty) \to (0,+\infty)$ non-decreasing |
| $g(v)$ | $g:\mathbb{R} \to \mathbb{R}$, $g'(v) > 0$ for $v \in \mathbb{R}$ | | $g(v) \equiv v$ | $g(v) \equiv v$ |
| $V(r)$ | $\equiv v^*$ | | $= v^* \sigma(r)$ | $\equiv 0$ |
| $S(v,r)$ | $\equiv 1$ | | $= -\dfrac{v_{\max}^2 v^* \sigma'(r)}{(v_{\max} - v) v}$ | $\equiv 1$ |
| $f(v)$ | $f:(-v^*, v_{\max} - v^*) \to \mathbb{R}$ | $= v$ | $= v$ | $= v$ |
| $Q(v,r)$ | $= v_{\max}^2 \dfrac{v_{\max} v + (v_{\max} - 2v) v^*}{2(v_{\max} - v)^2 v^2}$ | $\equiv 1$ | $= v_{\max}^2 \dfrac{v_{\max} v + (v_{\max} - 2v) v^* \sigma(r)}{2(v_{\max} - v)^2 v^2}$ | $\equiv 1$ |
| $J$ | $\equiv 1$ | $= \gamma - \dfrac{G}{v^*} + \dfrac{v_{\max} \ell(G)}{v^*(v_{\max} - v^*)}$ | $\equiv \gamma$ | $= \kappa(\rho, v)$ |
| Comments | $v f(v) > 0$ for $v \neq 0$ | $\ell: \mathbb{R} \to \mathbb{R}$, $\max(0,x) \leq \ell(x)$ for $x \in \mathbb{R}$ | $\sigma: \mathbb{R} \to (0,1]$ non-increasing with $\sigma(x) = 1$ for $x \leq \varepsilon$ and $\sigma(x) x \leq M$ for $x \geq \varepsilon$, $\|\sigma'\|_\infty + \|\sigma''\|_\infty < +\infty$ | $\kappa:(0,+\infty) \times \mathbb{R} \to \mathbb{R}_+$ |